\documentclass[aop,MSNbibl,seceqn,dvips]{arximspdf}

\doi{10.1214/12-AOP780} 
\volume{41}
\issue{4}
\pubyear{2013}
\firstpage{2479}
\lastpage{2512}

\makeatletter
\setattribute{abstract}    {skip} {24\p@}
\setattribute{frontmatter} {skip} {0\p@ plus 3\p@ minus 3\p@}

\newtheorem{Theorem}{Theorem}[section]
\newtheorem{Corollary}[Theorem]{Corollary}
\newtheorem{Proposition}[Theorem]{Proposition}
\newtheorem{Lemma}[Theorem]{Lemma}

\newproclaim{Remark}[Theorem]{Remark}

\newcommand{\E}{{\mathbf E}}
\newcommand{\R}{{\mathbf R}}
\renewcommand{\P}{{\mathbf P}}

\newcommand{\Var}{\operatorname{Var}}

\newcommand{\ep}{\varepsilon}

\makeatother

\begin{document}
\begin{frontmatter}

\title{Rate of convergence and Edgeworth-type expansion in the entropic
central limit theorem\thanksref{T1}}
\runtitle{Entropic central limit theorem}

\begin{aug}
\author[A]{\fnms{Sergey G.} \snm{Bobkov}\ead[label=e1]{bobkov@math.umn.edu}},
\author[B]{\fnms{Gennadiy P.} \snm{Chistyakov}\ead[label=e2]{chistyak@math.uni-bielefeld.de}}\\
\and
\author[B]{\fnms{Friedrich} \snm{G\"otze}\corref{}\ead[label=e3]{goetze@mathematik.uni-bielefeld.de}}
\runauthor{S. G. Bobkov, G. P. Chistyakov and F. G\"otze}
\affiliation{University of Minnesota, University of Bielefeld and
University of Bielefeld}
\address[A]{S. G. Bobkov\\
School of Mathematics\\
University of Minnesota\\
127 Vincent Hall, 206 Church St. S.E.\\
Minneapolis, Minnesota 55455\\
USA\\
\printead{e1}} 
\address[B]{G. P. Chistyakov\\
F. G\"otze\\
Fakult\"at f\"ur Mathematik\\
Universit\"at Bielefeld\\
Postfach 100131\\
33501 Bielefeld\\
Germany\\
\printead{e2}\\
\hphantom{E-mail: }\printead*{e3}}
\end{aug}

\thankstext{T1}{Supported in part by NSF Grant DMS-11-06530 and SFB 701.}

\received{\smonth{4} \syear{2011}}
\revised{\smonth{5} \syear{2012}}

%
\begin{abstract}
An Edgeworth-type expansion is established for the entropy distance to
the class of normal distributions of sums of i.i.d. random variables or
vectors, satisfying minimal moment conditions.
\end{abstract}

%

\end{frontmatter}

\section{Introduction}\label{sec1}

Let $(X_n)_{n \geq1}$ be independent, identically distributed
random variables with mean $\E X_1 = 0$ and variance $\Var(X_1) = 1$.
According to the central limit theorem, the normalized sums
\[
Z_n = \frac{X_1 + \cdots+ X_n}{\sqrt{n}}
\]
are weakly convergent in distribution to the standard normal law
$Z_n \Rightarrow Z$, where $Z \sim N(0,1)$ with density
$\varphi(x) = \frac{1}{\sqrt{2 \pi}} e^{-x^2/2}$.
A much\vspace*{1pt} stronger statement (when applicable)---the entropic central limit
theorem---states that, if for some $n_0$, or equivalently, for all
$n \geq n_0$, the random variables $Z_n$ have absolutely continuous
distributions with finite
entropies $h(Z_n)$, then these entropies converge,
%
%
\begin{equation}\label{equ1.1}
h(Z_n) \rightarrow h(Z)\qquad\mbox{as } n \rightarrow\infty.
\end{equation}
This theorem is due to Barron~\cite{Ba}. Some weaker variants of the theorem
in case of regularized distributions were known before; they go back to
the work of Linnik~\cite{L}, initiating an information-theoretic approach
to the central limit theorem.

To clarify in which sense (\ref{equ1.1}) is strong, recall that,
if a random variable $X$ with finite second moment has a density $p(x)$,
its entropy
\[
h(X) = -\int_{-\infty}^{+\infty} p(x) \log p(x) \,dx
\]
is well defined and is bounded from above by the entropy of the normal random
variable $Z$, having the same mean $a$ and the same variance $\sigma^2
$ as~$X$.
Note that the value $h(X) = -\infty$ is possible. The relative entropy
\[
D(X) = D(X\|Z) = h(Z) - h(X) = \int_{-\infty}^{+\infty}
p(x) \log\frac{p(x)}{\varphi_{a,\sigma}(x)} \,dx,
\]
where $\varphi_{a,\sigma}$ stands for the density of $Z$, is
nonnegative and
serves as kind of a distance to the class of normal laws, or to Gaussianity.
This quantity does
not depend on the mean or the variance of $X$, and can be related to
the total
variation distance between the distributions of $X$ and $Z$ by virtue
of the
Pinsker-type inequality $D(X) \geq\frac{1}{2} \|F_X - F_Z\|_{\mathrm{TV}}^2$.
This already shows that the entropic convergence (\ref{equ1.1}) is
stronger than
convergence in the total variation norm.

Thus, the entropic central limit theorem may be reformulated as
$D(Z_n) \rightarrow0$, as long as $D(Z_{n_0}) < +\infty$ for some $n_0$.
This property itself gives rise to a number of intriguing questions,
such as to the type and the rate of convergence. In particular, it has been
proved only recently that the sequence $h(Z_n)$ is nondecreasing,
so that $D(Z_n) \downarrow0$; cf.~\cite{A-B-B-N1,B-M}.
This leads to the question as to the precise rate of $D(Z_n)$ tending
to zero;
however, not much seems to be known about this problem. The best
results in
this direction are due to Artstein et al.~\cite{A-B-B-N2} and
to Barron and Johnson~\cite{B-J}. In the i.i.d. case as above, these authors
have obtained
an expected asymptotic bound $D(Z_n) = O(1/n)$ under the hypothesis
that the
distribution of $X_1$ admits an analytic inequality of Poincar\'e-type
(in~\cite{B-J}, a restricted Poincar\'e inequality is used).
These inequalities involve a large variety of ``nice'' probability
distributions which necessarily have a finite exponential moment.

The aim of this paper is to study the rate of $D(Z_n)$, using moment conditions
$\E|X_1|^s < +\infty$ with fixed values $s \geq2$, which are
comparable to those
required for classical Edgeworth-type approximations in the Kolmogorov distance.
The cumulants
\[
\gamma_r = i^{-r} \,\frac{d^r}{dt^r} \log\E
e^{itX_1}\bigg|_{t=0}
\]
are then
well defined for all $r \leq[s]$ (the integer part of $s$), and one may
introduce the functions
%
%
\begin{equation}\label{equ1.2}
q_k(x) = \varphi(x) \sum H_{k + 2j}(x)
\frac{1}{r_1!\cdots r_k!} \biggl(\frac{\gamma_3}{3!} \biggr)^{r_1}
\cdots\biggl(
\frac{\gamma_{k+2}}{(k+2)!} \biggr)^{r_k}
\end{equation}
involving the Chebyshev--Hermite polynomials $H_k$. The summation in
(\ref{equ1.2}) runs
over all nonnegative integer solutions $(r_1,\ldots,r_k)$ to the equation
$r_1 + 2 r_2 + \cdots+ k r_k = k$, and one uses the notation
$j = r_1 + \cdots+ r_k$.

The functions $q_k$ are defined for $k = 1,\ldots,[s]-2$. They appear in
Edgeworth-type expansions including the local limit theorem, where
$q_k$ are used to construct the approximation of the densities of
$Z_n$. These results can be applied to obtain an expansion in powers of
$1/n$ for the distance $D(Z_n)$. For a multidimensional version of the
following Theorem~\ref{Theorem1.1} for moments of integer order $s
\geq2$, see Theorem~\ref{Theorem6.1} below.

%
\begin{Theorem}\label{Theorem1.1} Let $\E|X_1|^s < +\infty$ $(s \geq
2)$, and assume
$D(Z_{n_0}) < +\infty$, for some $n_0$. Then
%
%
\begin{equation}\label{equ1.3}
D(Z_n)= \frac{c_1}{n} + \frac{c_2}{n^2} + \cdots+
\frac{c_{[(s-2)/2]}}{n^{[(s-2)/2]}} + o \bigl((n \log n)^{-(s-2)/2}
\bigr).
\end{equation}
Here
%
%
\begin{equation}\label{equ1.4}
c_j = \sum_{k=2}^{2j}
\frac{(-1)^k}{k(k-1)} \sum\int_{-\infty}^{+\infty}
q_{r_1}(x) \cdots q_{r_k}(x) \,\frac{dx}{\varphi(x)^{k-1}},
\end{equation}
where the summation runs over all positive integers
$(r_1,\ldots,r_k)$ such that $r_1 + \cdots+ r_k = 2j$.
\end{Theorem}


Each coefficient $c_j$ in (\ref{equ1.3}) represents a certain
polynomial in the cumulants
$\gamma_3,\ldots,\gamma_{2j+1}$. For example, $c_1 = \frac{1}{12}
\gamma_3^2$,
and in the case $s=4$, (\ref{equ1.3}) gives
%
%
\begin{equation}\label{equ1.5}
D(Z_n) = \frac{1}{12 n} \bigl(\E X_1^3
\bigr)^2 + o \biggl(\frac{1}{n \log n} \biggr) \qquad\bigl( \E
X_1^4 < +\infty\bigr).
\end{equation}
Thus, under the 4th moment condition,
we have $D(Z_n) \leq\frac{C}{n}$, where the constant depends on the underlying
distribution. This has been conjectured by Johnson~\cite{J}, page 49.
Actually, the constant
$C$ may be expressed in terms of $\E X_1^4$ and $D(X_1)$, only.

When $s$ varies in the range $4 \leq s \leq6$, the leading linear term
in (\ref{equ1.5})
will be unchanged, while the remainder term improves and satisfies
$O(\frac{1}{n^2})$ in case $\E X_1^6 < +\infty$.
But for $s=6$, the result involves the subsequent coefficient $c_2$
which depends
on $\gamma_3,\gamma_4$ and $\gamma_5$.
In particular, if $\gamma_3 = 0$, we have $c_2 = \frac{1}{48} \gamma
_4^2$, thus
\[
D(Z_n) = \frac{1}{48 n^2} \bigl(\E X_1^4
- 3 \bigr)^2 + o \biggl(\frac{1}{(n \log n)^2} \biggr)\qquad\bigl
( \E
X_1^3 = 0, \E X_1^6 < +
\infty\bigr).
\]
More generally, representation (\ref{equ1.3}) simplifies if the first $k-1$
moments of $X_1$ coincide with the corresponding moments of $Z \sim N(0,1)$.

%
\begin{Corollary}\label{Corollary1.2} Let $\E|X_1|^s < +\infty$ $(s \geq
4)$, and
assume that
$D(Z_{n_0}) < +\infty$, for some $n_0$. Given $k = 3,4,\ldots,[s]$,
assume that
$\gamma_j = 0$ for all $3 \leq j < k$. Then
%
%
\begin{equation}\label{equ1.6}
D(Z_n) = \frac{\gamma_k^2}{2 k!}\cdot\frac{1}{n^{k-2}} + O \biggl(
\frac{1}{n^{k-1}} \biggr) + o \biggl(\frac{1}{(n \log n)^{(s-2)/2}}
\biggr).
\end{equation}
\end{Corollary}

Johnson had noticed (though in terms of the standardized Fisher
information, see
\cite{J}, Lemma 2.12) that if $\gamma_k \neq0$, $D(Z_n)$ cannot be of
smaller order than
$n^{-(k-2)}$.

Note that when $\E X_1^{2k} < +\infty$, the $o$-term may be removed in the
representation~(\ref{equ1.6}). On the other hand, when $k > \frac
{s+2}{2}$, the $o$-term
will dominate the $n^{-(k-2)}$-term, and we can only conclude
that $D(Z_n) = o ((n \log n)^{-(s-2)/2} )$.

As for the missing range $2 \leq s < 4$, here there are no coefficients $c_j$
appearing in the sum (\ref{equ1.3}), and Theorem~\ref{Theorem1.1} just
tells us that
%
%
\begin{equation}\label{equ1.7}
D(Z_n) = o \biggl(\frac{1}{(n \log n)^{(s-2)/2}} \biggr).
\end{equation}
This bound is worse than the rate $1/n$. In particular, it only gives
$D(Z_n) = o(1)$
for $s=2$, which is the statement of Barron's theorem. In fact, in this case
the entropic distance to normality may decay to zero at an arbitrarily
slow rate.
In case of a finite 3rd absolute moment,
$D(Z_n) = o(\frac{1}{\sqrt{n \log n}})$.
To see that this and that the more general relation (\ref{equ1.7})
cannot be improved
with respect to the powers of $1/n$, we prove:

%
\begin{Theorem}\label{Theorem1.3} Let $\eta>1$. Given $2 < s < 4$,
there exists
a sequence
of independent, identically distributed random variables $(X_n)_{n \geq
1}$ with
$\E|X_1|^s < +\infty$, such that $D(X_1) < +\infty$ and
\[
D(Z_n) \geq\frac{c}{(n \log n)^{(s-2)/2} (\log n)^{\eta}},\qquad n \geq n_1(X_1),
\]
with a constant $c=c(\eta,s)>0$, depending on $\eta$ and $s$, only.
\end{Theorem}

Known bounds on the entropy are commonly based on Bruijn's identity
which may be used to represent the entropic distance to normality as an
integral of the Fisher information for regularized distributions; cf.
\cite{Ba}. However, it is not clear how to reach exact asymptotics with
this approach. The proofs of Theorems~\ref{Theorem1.1} and
\ref{Theorem1.3} stated above rely upon classical tools and results in
the theory of sums of independent summands including Edgeworth-type
expansions for convolution of densities formulated as local limit
theorems with nonuniform remainder bounds. For noninteger values
of $s$, the authors had to complete the otherwise extensive literature
by recent, technically rather involved results based on fractional
differential calculus; see~\cite{B-C-G1,B-C-G2}. Our approach applies to random
variables in higher dimension as well and to nonidentical distributions
for summands with uniformly bounded $s$th moments.

We start with the description of a truncation-of-density argument,
which allows us to reduce many questions about bounding the entropic
distance to the case of bounded densities (Section~\ref{sec2}). In
Section~\ref{sec3} we discuss known results about Edgeworth-type
expansions that will be used in the proof of Theorem~\ref{Theorem1.1}.
Main steps of the proofs are based on it in Sections~\ref{sec4} and~\ref{sec5}.
All auxiliary results cover the scheme of i.i.d. random
vectors in $\R^d$ as well (however, with integer values of~$s$) and are
finalized in Section~\ref{sec6} to obtain multidimensional variants of
Theorem~\ref{Theorem1.1} and Corollary~\ref{Corollary1.2}.
Sections~\ref{sec7} and~\ref{sec8} are devoted to lower bounds on the entropic
distance to normality for a special class of probability distributions
on the real line that are used in the proof Theorem~\ref{Theorem1.3}.

\section{Binomial decomposition of convolutions}\label{sec2}

First let us comment on the assumptions in Theorem~\ref{Theorem1.1}.
It may happen that $X_1$ has a singular distribution, but the distribution
of $X_1 + X_2$ and of all next sums $S_n = X_1 + \cdots+ X_n$ ($n \geq2$)
are absolutely continuous; cf.~\cite{T}.

If it exists, the density $p$ of $X_1$ may or may not be bounded. In
the first case,
all the entropies $h(S_n)$ are finite. If $p$ is unbounded, it may
happen that
all $h(S_n)$ are infinite, even if $p$ is compactly supported.
But if $h(S_n)$ is finite for some $n = n_0$ then, for all $n \geq n_0$,
entropies are finite; see~\cite{Ba} for specific examples.

Denote by $p_n(x)$ the density of $Z_n = S_n/\sqrt{n}$ (when it exists).
Since it is desirable to work with bounded densities, we will slightly modify
$p_n$ at the expense of a small change in the entropy. Variants of the next
construction are well known; see, for example,~\cite{S-M,I-L},
where the
central limit theorem
was studied with respect to the total variation distance.
Without any extra efforts, we may assume that $X_n$ take values in $\R^d$
which we equip with the usual inner product $\langle\cdot,\cdot
\rangle
$ and
the Euclidean norm $|\cdot|$. For simplicity, we describe the construction
in the situation, where $X_1$ has a density $p(x)$; cf. Remark \ref
{Remark2.5} on
appropriate modifications in the general case.

Let $m_0 \geq0$ be a fixed integer. (For the purposes of Theorem \ref
{Theorem1.1},
one may take $m_0 = [s] + 1$.)

If $p$ is bounded, we put $\widetilde p_n(x) = p_n(x)$ for all $n \geq1$.
Otherwise, the integral
%
%
\begin{equation}\label{equ2.1}
b = \int_{p(x) > M} p(x) \,dx
\end{equation}
is positive for all $M>0$.
Choose $M$ to be sufficiently large to satisfy, for example, $0 < b <
\frac{1}{2}$; cf. Remark~\ref{Remark2.4}.
In this case (when $p$ is unbounded), consider the decomposition
%
%
\begin{equation}\label{equ2.2}
p(x) = (1-b) \rho_1(x) + b \rho_2(x),
\end{equation}
where $\rho_1$, $\rho_2$ are the normalized restrictions of $p$
to the sets $\{p(x) \leq M\}$ and $\{p(x) > M\}$, respectively. Hence,
for the convolutions we have a binomial decomposition
\[
p^{*n} = \sum_{k=0}^n
C_n^k (1-b)^k b^{n-k}
\rho_1^{*k} * \rho_2^{*(n-k)}.
\]
For $n \geq m_0 + 1$, we split the above sum into the two parts, so that
$p^{*n} = \rho_{n1} + \rho_{n2}$ with
\begin{eqnarray*}
\rho_{n1} &=& \sum_{k = m_0 + 1}^n
C_n^k (1-b)^k b^{n-k}
\rho_1^{*k} * \rho_2^{*(n-k)},\\
\rho_{n2} &=& \sum_{k=0}^{m_0}
C_n^k (1-b)^k b^{n-k}
\rho_1^{*k} * \rho_2^{*(n-k)}.
\end{eqnarray*}
Note that, whenever $b < b_1 < \frac{1}{2}$,
%
%
\begin{eqnarray}\label{equ2.3}
\ep_n &\equiv&\int\rho_{n2}(x) \,dx = \sum
_{k=0}^{m_0} C_n^k
(1-b)^k b^{n-k} \nonumber\\[-8pt]\\[-8pt]
&\leq& n^{m_0} b^{n-m_0} = o
\bigl(b_1^n\bigr)\qquad\mbox{as } n \rightarrow\infty.\nonumber
\end{eqnarray}
Finally define
%
%
\begin{equation}\label{equ2.4}
\widetilde p_n(x) = p_{n1}(x) = \frac{1}{1 - \ep_n}
n^{d/2} \rho_{n1} (x\sqrt{n} )
\end{equation}
and similarly $p_{n2}(x) = \frac{1}{\ep_n} n^{d/2} \rho_{n2}(x\sqrt
{n})$. Thus,
we have the desired decomposition
%
%
\begin{equation}\label{equ2.5}
p_n(x) = (1 - \ep_n) p_{n1}(x) +
\ep_n p_{n2}(x).
\end{equation}

The probability densities $p_{n1}(x)$ are bounded and
provide an approximation for $p_n(x) = n^{d/2} p^{*n}(x\sqrt{n})$ in
total variation.
In particular, from (\ref{equ2.3})--(\ref{equ2.5}) it follows that
\[
\int\bigl|p_{n1}(x) - p_n(x)\bigr| \,dx < 2^{-n}
\]
for all $n$ large enough. One of the immediate consequences of this estimate
is the bound
%
%
\begin{equation}\label{equ2.6}
\bigl|v_{n1}(t) - v_n(t)\bigr| < 2^{-n} \qquad\bigl(t \in
\R^d\bigr)
\end{equation}
for the characteristic functions
$v_n(t) = \int e^{i\langle t,x\rangle} p_n(x) \,dx$ and
$v_{n1}(t) = \int e^{i\langle t,x\rangle} \*p_{n1}(x) \,dx$,
corresponding to the densities $p_n$ and $p_{n1}$.

This property may be sharpened in case of finite moments.

%
\begin{Lemma}\label{Lemma2.1} If $\E|X_1|^s < +\infty$ $(s \geq0)$, then
for all $n$ large enough,
\[
\int\bigl(1 + |x|^s\bigr) \bigl|\widetilde p_n(x) -
p_n(x)\bigr| \,dx < 2^{-n}.
\]
In particular, (\ref{equ2.6}) also holds for all partial derivatives of
$v_{n1}$ and $v_n$ up to order $m = [s]$.
\end{Lemma}

\begin{pf} By definition (\ref{equ2.5}),
$|p_{n1}(x) - p_n(x)| \leq\ep_n (p_{n1}(x) + p_{n2}(x))$, hence
\begin{eqnarray*}
\int|x|^s \bigl|p_{n1}(x) - p_n(x)\bigr| \,dx &\leq&
\frac{\ep_n}{1 - \ep_n} n^{-s/2} \int|x|^s \rho_{n1}(x)
\,dx\\
&&{} + n^{-s/2} \int|x|^s \rho_{n2}(x) \,dx.
\end{eqnarray*}

Let $U_1,U_2,\ldots$ be independent copies of $U$ and $V_1,V_2,\ldots$
be independent copies of $V$ (that are also independent of $U_n$'s), where
$U$ and $V$ are random vectors with densities $\rho_1$ and $\rho_2$,
respectively. From (\ref{equ2.2})
\[
\beta_s \equiv\E|X_1|^s = (1-b) \E
|U|^s + b \E|V|^s,
\]
so $\E|U|^s \leq\beta_s/b$ and $\E|V|^s \leq\beta_s/b$ (using $b
< \frac{1}{2}$).
Therefore, for the normalized sums
\[
R_{k,n} = \frac{1}{\sqrt{n}} (U_1 + \cdots+
U_k + V_1 + \cdots+ V_{n-k}),\qquad 0 \leq k \leq n,
\]
we have $\E|R_{k,n}|^s \leq\frac{\beta_s}{b} n^{s/2}$, if $s \geq
1$, and
$\E|R_{k,n}|^s \leq\frac{\beta_s}{b} n^{1 - (s/2)}$, if $0 \leq s
\leq1$.
Hence, by the definition of $\rho_{n1}$ and $\rho_{n2}$,
\begin{eqnarray*}
\int|x|^s \rho_{n1}(x) \,dx &=& n^{s/2} \sum
_{k=m_0+1}^{n} C_n^k
(1-b)^k b^{n-k} \E|R_{k,n}|^s \leq
\frac{\beta_s}{b} n^{s+1},
\\
\int|x|^s \rho_{n2}(x) \,dx &=& n^{s/2} \sum
_{k=0}^{m_0} C_n^k
(1-b)^k b^{n-k} \E|R_{k,n}|^s \leq
\frac{\beta_s}{b} n^{s+1} \ep_n.
\end{eqnarray*}
It remains to apply estimate (\ref{equ2.3}) on $\ep_n$, and Lemma \ref
{Lemma2.1} follows.
\end{pf}

We need to extend the assertion of Lemma~\ref{Lemma2.1} to the relative
entropies
with respect to the standard normal distribution on $\R^d$ with density
$\varphi(x) = (2\pi)^{-d/2} e^{-|x|^2/2}$. Thus put
\[
D_n = \int p_n(x) \log\frac{p_n(x)}{\varphi(x)} \,dx,\qquad \widetilde
D_n = \int\widetilde p_n(x) \log\frac{\widetilde
p_n(x)}{\varphi(x)} \,dx.
\]

%
\begin{Lemma}\label{Lemma2.2} If $X_1$ has a finite second moment and
finite entropy,
then $|\widetilde D_n - D_n| < 2^{-n}$, for all $n$ large enough.
\end{Lemma}

First, we collect a few elementary properties of the convex function
$L(u) = u \log u$ ($u \geq0$).

%
\begin{Lemma}\label{Lemma2.3}
For all $u,v \geq0$ and $0 \leq\ep\leq1$:

\begin{longlist}[(a)]
\item[(a)]
$L((1 - \ep) u + \ep v) \leq(1-\ep) L(u) + \ep L(v)$;

\item[(b)] $L((1 - \ep) u + \ep v) \geq(1-\ep) L(u) + \ep L(v) +
u L(1-\ep) + v L(\ep)$;

\item[(c)] $L((1 - \ep) u + \ep v) \geq(1-\ep) L(u) - \frac{1}{e} u -
\frac{1}{e}$.
\end{longlist}
\end{Lemma}

The first assertion is just Jensen's inequality applied to $L$.
By the convexity of~$L$, for each $y \geq0$, the function
$L(x+y)-L(x)$ is
increasing in $x \geq0$. Hence, $L(x+y)-L(x) \geq L(y)$, which is (b)
for $x = (1 - \ep) u$ and $y = \ep v$. Similarly, using
$L \geq- \frac{1}{e}$, we obtain (c).

\begin{pf*}{Proof of Lemma~\ref{Lemma2.2}}
Assuming that $p$ is (essentially) unbounded, define
\[
D_{nj} = \int p_{nj}(x) \log\frac{p_{nj}(x)}{\varphi(x)} \,dx\qquad
(j=1,2),
\]
so that $\widetilde D_n = D_{n,1}$. By Lemma~\ref{Lemma2.3}(a),
$D_n \leq(1 - \ep_n)D_{n1} + \ep_n D_{n2}$. On the other hand, by (b),
\[
D_n \geq\bigl((1 - \ep_n)D_{n1} +
\ep_n D_{n2} \bigr) + \ep_n \log
\ep_n + (1-\ep_n) \log(1-\ep_n).
\]
In view of (\ref{equ2.3}), the two estimates give
%
%
\begin{equation}\label{equ2.7}
|D_{n1} - D_n| < C(n + D_{n1} +
D_{n2}) b_1^n,
\end{equation}
which holds for all $n \geq1$ with some constant $C$. In addition, by
the inequality
in (c) with $\ep= b$, from (\ref{equ2.2}) it follows that
%
%
\begin{equation}\label{equ2.8}\qquad
D(X_1\|Z) = \int L \biggl(\frac{p(x)}{\varphi(x)} \biggr) \varphi
(x) \,dx
\geq(1-b) \int\rho_1(x) \log\frac{\rho_1(x)}{\varphi(x)} \,dx -
\frac{2}{e},
\end{equation}
where $Z$ denotes a standard normal random vector in $\R^d$.
By the same reasoning,
%
%
\begin{equation}\label{equ2.9}
D(X_1\|Z) \geq b \int\rho_2(x) \log\frac{\rho_2(x)}{\varphi(x)}
\,dx - \frac{2}{e}.
\end{equation}

Now, by the convexity of the function $L(u) = u\log u$,
\begin{eqnarray*}
D_{n1} & \leq& \frac{1}{1-\ep_n} \sum_{k=m_0+1}^{n}
C_n^k (1-b)^k b^{n-k} \int
r_{k,n}(x) \log\frac{r_{k,n}(x)}{\varphi(x)} \,dx,
\\
D_{n2} & \leq& \frac{1}{\ep_n} \sum_{k=0}^{m_0}
C_n^k (1-b)^k b^{n-k} \int
r_{k,n}(x) \log\frac{r_{k,n}(x)}{\varphi(x)} \,dx,
\end{eqnarray*}
where $r_{k,n}$ are densities of the normalized sums $R_{k,n}$
from the proof of Lem\-ma~\ref{Lemma2.1}. Here each integral may also be
written as
%
%
\begin{equation}\label{equ2.10}\qquad
\int r_{k,n}(x) \log\frac{r_{k,n}(x)}{\varphi(x)} \,dx = \int L
\bigl(r_{k,n}(x)\bigr) \,dx + \frac{d}{2} \log(2\pi) +
\frac{1}{2} \E|R_{k,n}|^2.
\end{equation}
We have $\E|R_{k,n}|^2 \leq\frac{\beta^2}{b} n$, as noticed in the
proof of Lemma~\ref{Lemma2.1}. In addition, by the convexity of $L$,
there is a general
inequality
\[
\int L\bigl((f * g) (x)\bigr) \,dx \leq\int L\bigl(f(x)\bigr) \,dx
\]
valid for the convolution of any two probability densities $f$ and $g$
on $\R^d$
(if the integrals exist). In particular,
\[
\int L\bigl(r_{k,n}(x)\bigr) \,dx \leq\frac{d}{2} \log n + \max
\biggl\{\int L\bigl(\rho_1(x)\bigr) \,dx, \int L\bigl(
\rho_2(x)\bigr) \,dx \biggr\},
\]
which may actually be sharpened in case $1<k<n$ by replacing $\max$
with $\min$.
By (\ref{equ2.8}) and (\ref{equ2.9}), the integrals on the right-hand
side are finite, thus
the integrals on
the left-hand side of (\ref{equ2.10}) are bounded by $Cn$ with some
constant $C$.
Hence, a similar bound also holds for $D_{nj}$,
and it remains to apply (\ref{equ2.7}). Lemma~\ref{Lemma2.2} is proved.
\end{pf*}

%
\begin{Remark}\label{Remark2.4} If $X_1$ has a finite second moment and
$D(X_1) <
+\infty$,
the truncation level $M$ in (\ref{equ2.1}) can be chosen explicitly in
terms of
$b$ using
the entropic distance $D(X_1)$ and $\sigma^2 = \operatorname{det}(\Sigma)$,
where $\Sigma$ is the covariance matrix of $X_1$.

Indeed, putting $a = \E X_1$ and using an elementary inequality
$t\log(1+t) \leq t \log t + 1$ ($t \geq0$), we have an upper estimate
\begin{eqnarray*}
\int p \log\biggl(1 + \frac{p}{\varphi_{a,\Sigma}} \biggr) \,dx & = & \int
\frac{p}{\varphi_{a,\Sigma}} \log\biggl(1 + \frac{p}{\varphi
_{a,\Sigma}} \biggr)
\varphi_{a,\Sigma} \,dx
\\
& \leq& \int p \log\frac{p}{\varphi_{a,\Sigma}} \,dx + 1 = D(X_1) + 1.
\end{eqnarray*}
On the other hand, the original expression majorizes
\[
\int_{\{p(x) > M\}} p(x) \log\frac{M}{\varphi_{a,\Sigma}(x)} \,dx
\geq b \log
\bigl(M\sigma(2\pi)^{d/2} \bigr),
\]
hence
\[
M \leq\frac{1}{\sigma(2\pi)^{d/2}} e^{(D(X_1) + 1)/b}.
\]
\end{Remark}

%
\begin{Remark}\label{Remark2.5}
If $Z_n$ have absolutely continuous distributions
with finite
entropies for $n \geq n_0 > 1$, the above construction should be
properly modified.\vadjust{\goodbreak}

Namely, one may put $\widetilde p_n = p_n$, if $p_n$ are bounded, and otherwise
apply the same decomposition (\ref{equ2.2}) to $p_{n_0}$ in place of
$p$. As a result,
for any $n = An_0 + B$ ($A \geq1$, $0 \leq B \leq n_0 - 1$),
the partial sum $S_n$ will have the density
\[
r_n(x) = \sum_{k=0}^A
C_A^k (1-b)^k b^{A-k} \int\bigl(
\rho_1^{*k} * \rho_2^{*(A-k)} \bigr) (x
- y) \,dF_B(y),
\]
where $F_B$ is the distribution of $S_B$. For $A \geq m_0 + 1$, split
the above sum
into the two parts with summation over $m_0 + 1 \leq k \leq A$ and
$0 \leq k \leq m_0$, respectively, so that $r_n = \rho_{n1} + \rho_{n2}$.
Then, like in (\ref{equ2.4}) and for the same sequence $\ep_n$
described in~(\ref{equ2.3}), define
\[
\widetilde p_n(x) = \frac{1}{1 - \ep_n} n^{d/2}
\rho_{n1} (x\sqrt{n} ).
\]
Clearly, these densities are bounded and approximate $p_n(x)$ in total
variation.
In particular, for all sufficiently large $n$, they satisfy the
estimates that are
similar to the estimates in Lemmas~\ref{Lemma2.1} and~\ref{Lemma2.2}.
\end{Remark}

\section{Edgeworth-type expansions}\label{sec3}

Let $(X_n)_{n \geq1}$ be independent, identically distributed random
variables with mean $\E X_1 = 0$ and variance \mbox{$\Var(X_1) = 1$}.
In this section we collect some auxiliary results about Edgeworth-type
expansions
both for the distribution functions $F_n(x) = \P\{Z_n \leq x\}$ and the
densities
$p_n(x)$ of the normalized sums $Z_n = S_n/\sqrt{n}$, where $S_n = X_1
+ \cdots+ X_n$.

If the absolute moment $\E|X_1|^s$ is finite for a given $s \geq2$
and $m = [s]$,
define
%
%
\begin{equation}\label{equ3.1}
\varphi_m(x) = \varphi(x) + \sum_{k=1}^{m-2}
q_k(x) n^{-k/2}
\end{equation}
with the functions $q_k$ described in (\ref{equ1.2}). Introduce as well
%
%
\begin{equation}\label{equ3.2}
\Phi_m(x) = \int_{-\infty}^x
\varphi_m(y) \,dy = \Phi(x) + \sum_{k=1}^{m-2}
Q_k(x) n^{-k/2}.
\end{equation}
Similar to (\ref{equ1.2}), the functions $Q_k$ have an explicit
description involving
the cumulants $\gamma_3,\ldots,\gamma_{k+2}$ of $X_1$. Namely,
\[
Q_k(x) = -\varphi(x) \sum H_{k + 2j-1}(x)
\frac{1}{r_1!\cdots r_k!} \biggl(\frac{\gamma_3}{3!} \biggr)^{r_1}
\cdots\biggl(
\frac{\gamma_{k+2}}{(k+2)!} \biggr)^{r_k},
\]
where the summation is carried out over all nonnegative integer
solutions $(r_1,\ldots,r_k)$ to the equation
$r_1 + 2 r_2 + \cdots+ k r_k = k$ with $j = r_1 + \cdots+ r_k$; cf.,
for example,~\cite{B-RR} or~\cite{Pe2} for details.

%
\begin{Theorem}\label{Theorem3.1} Assume that
$\limsup_{|t| \rightarrow+\infty} |\E e^{it X_1}| < 1$.
If $\E|X_1|^s < +\infty$ $(s \geq2)$, then as $n \rightarrow\infty$,
uniformly for all $x$,
%
%
\begin{equation}\label{equ3.3}
\bigl(1 +|x|^s\bigr) \bigl(F_n(x) -
\Phi_m(x) \bigr) = o \bigl(n^{-(s-2)/2} \bigr).
\end{equation}
\end{Theorem}


For $2 \leq s < 3$ and $m=2$, there are no expansion terms in the sum
(\ref{equ3.2}), and hence
$\Phi_2(x) = \Phi(x)$ is the distribution function
of the standard normal law. In this case, (\ref{equ3.3}) becomes
%
%
\begin{equation}\label{equ3.4}
\bigl(1 +|x|^s\bigr) \bigl(F_n(x) - \Phi(x) \bigr) = o
\bigl(n^{-(s-2)/2} \bigr).
\end{equation}
In fact, in this case Cramer's condition on the characteristic function
of $X_1$
is not used. The result was obtained by Osipov and Petrov~\cite{O-P}; cf.
also~\cite{Bi} where (\ref{equ3.4}) is established with $O$.

In the case $s \geq3$ Theorem~\ref{Theorem3.1} can be found in~\cite{Pe2}
(Theorem 2, Chapter~VI, page 168). Note that when $s=m$ is
integer, relation (\ref{equ3.3}) without the factor $1 +|x|^m$ represents
the classical Edgeworth expansion. It is essentially due to Cram\'er and
is described in many papers and textbooks; cf.~\cite{E,F}.
However, the case of fractional values of $s$ is more delicate, especially
in the following local limit theorem.

%
\begin{Theorem}\label{Theorem3.2}
Let $\E|X_1|^s < +\infty$ $(s \geq2)$. Suppose $Z_{n_0}$ has a bounded
density for some $n_0$. Then for all sufficiently large $n$, the random
variables $Z_n$ have continuous bounded densities $p_n$ satisfying, as
$n \rightarrow\infty$,
%
%
\begin{equation}\label{equ3.5}
\bigl(1 + |x|^m\bigr) \bigl(p_n(x) -
\varphi_m(x) \bigr) = o \bigl(n^{-(s-2)/2} \bigr)
\end{equation}
uniformly for all $x$. Moreover,
%
%
\begin{eqnarray}\label{equ3.6}
&&
\bigl(1 + |x|^s\bigr) \bigl(p_n(x) -
\varphi_m(x) \bigr)\nonumber\\[-8pt]\\[-8pt]
&&\qquad= o \bigl(n^{-(s-2)/2} \bigr) +
\bigl(1+|x|^{s-m}\bigr) \bigl(O\bigl(n^{-(m-1)/2}\bigr)+ o
\bigl(n^{-(s-2)}\bigr) \bigr).\nonumber
\end{eqnarray}
\end{Theorem}

If $s = m$ is integer and $m \geq3$, Theorem~\ref{Theorem3.2} is well
known; then
(\ref{equ3.5}) and (\ref{equ3.6})
simplify to
%
%
\begin{equation}\label{equ3.7}
\bigl(1 + |x|^m\bigr) \bigl(p_n(x) -
\varphi_m(x) \bigr) = o \bigl(n^{-(m-2)/2} \bigr).
\end{equation}
In this formulation the result is due to Petrov~\cite{Pe1}; cf. \cite
{Pe2}, page
211, or~\cite{B-RR}, pa\-ge~192. Without the term $1 + |x|^m$,
relation (\ref{equ3.7}) goes back to the results of Cram\'er and
Gnedenko (cf.~\cite{G-K}).

In the general (fractional) case, Theorem~\ref{Theorem3.2} has recently been
obtained in~\cite{B-C-G1,B-C-G2}
by using the technique of Liouville fractional integrals and derivatives.
Assertion (\ref{equ3.6}) gives an improvement over (\ref{equ3.5}) on
relatively large intervals
of the real axis, and this is essential in the case of noninteger $s$.

An obvious weak point in Theorem~\ref{Theorem3.2} is that it requires the
boundedness of
the densities $p_n$, which is, however, necessary for conclusions, such as
(\ref{equ3.5}) or (\ref{equ3.7}). Nevertheless, this condition may be
removed, if we
replace $p_n$
by slightly modified densities $\widetilde p_n$.

%
\begin{Theorem}\label{Theorem3.3}
Let $\E|X_1|^s < +\infty$ $(s \geq2)$. Suppose
that, for all for all sufficiently large $n$, $Z_n$ have absolutely
continuous distributions
with densities $p_n$. Then there exist some bounded continuous densities
$\widetilde p_n$ such that:

\begin{longlist}[(a)]
\item[(a)]
the relations (\ref{equ3.5}) and (\ref{equ3.6}) hold true for
$\widetilde p_n$
instead of $p_n$;

\item[(b)] $\int_{-\infty}^{+\infty} (1 + |x|^s) |\widetilde p_n(x) -
p_n(x)| \,dx < 2^{-n}$,
for all sufficiently large $n$;

\item[(c)]
$\widetilde p_n(x) = p_n(x)$ almost everywhere, if $p_n$ is bounded
$(\mbox{a.e.})$.
\end{longlist}
\end{Theorem}

Here, property (c) is added to include Theorem~\ref{Theorem3.2} in
Theorem~\ref{Theorem3.3} as a particular case. Moreover, one can use
the densities $\widetilde p_n$ constructed in the previous section with
$m_0 = [s]+1$. We refer to~\cite{B-C-G1,B-C-G2} for detailed proofs.

This extended result allows us to immediately recover, for example,
the central limit theorem with respect to the total variation distance
(without the assumption of boundedness of $p_n$). Namely, we have
%
%
\begin{equation}\label{equ3.8}
\|F_n - \Phi_m\|_{\mathrm{TV}} = \int
_{-\infty}^{+\infty} \bigl|p_n(x) -
\varphi_m(x)\bigr| \,dx = o \bigl(n^{-(s-2)/2} \bigr).
\end{equation}
For $s=2$ and $\varphi_2(x) = \varphi(x)$, this statement corresponds
to a theorem of Prokhorov~\cite{Pr}, while for $s=3$ and
$\varphi_3(x) = \varphi(x) (1 + \gamma_3 \frac{x^3 - 3x}{6\sqrt
{n}} )$---to the result of Sirazhdinov and Mamatov~\cite{S-M}.

\section*{The multidimensional case}

Similar results are also available in the multidimensional case for
integer values
$s = m$. In the remaining part of this section,
let $(X_n)_{n \geq1}$ denote independent identically distributed random
vectors in the Euclidean space $\R^d$ with mean zero and identity
covariance matrix.

Assuming $\E|X_1|^m < +\infty$ for some integer $m \geq2$ (where now
$|\cdot|$
denotes the Euclidean norm), introduce
the cumulants $\gamma_\nu$ of $X_1$ and the associated cumulant
polynomials $\gamma_k(it)$ up to order $m$ by using the equality
\[
\frac{1}{k!} \,\frac{d^k} {du^k} \log\E e^{iu\langle t,X_1\rangle} \bigg|_{u=0}
= \frac{1}{k!} \gamma_k(it) = \sum
_{|\nu| = k} \gamma_\nu\frac{(it)^\nu}{\nu!} \qquad\bigl(k = 1,
\ldots, m, t \in\R^d \bigr).
\]
Here the summation runs over all $d$-tuples $\nu= (\nu_1,\ldots,\nu_d)$
with integer
components $\nu_j \geq0$ such that $|\nu| = \nu_1 + \cdots+ \nu_d =
k$. We also
write $\nu! = \nu_1! \cdots\nu_d!$ and use a standard notation for the
generalized
powers $z^\nu= z_1^{\nu_1} \cdots z_d^{\nu_d}$ of real or complex vectors
$z = (z_1,\ldots,z_d)$, which are treated as polynomials in $z$ of
degree $|\nu|$.

For $1 \leq k \leq m-2$, define the polynomials
%
%
\begin{equation}\label{equ3.9}
P_k(it) = \sum_{r_1 + 2 r_2 + \cdots+ k r_k = k}
\frac{1}{r_1!\cdots r_k!} \biggl(\frac{\gamma_3(it)}{3!} \biggr
)^{r_1} \cdots\biggl(
\frac{\gamma_{k+2}(it)}{(k+2)!} \biggr)^{r_k},
\end{equation}
where the summation is performed over all nonnegative integer
solutions
$(r_1,\ldots,r_k)$ to the equation $r_1 + 2 r_2 + \cdots+ k r_k = k$.

Furthermore, like in dimension one, define the approximating functions
$\varphi_m(x)$
on $\R^d$ by virtue of the equality (\ref{equ3.1}), where every $q_k$
is determined
by its Fourier transform
%
%
\begin{equation}\label{equ3.10}
\int e^{i\langle t,x\rangle} q_k(x) \,dx = P_k(it)
e^{-|t|^2/2}.
\end{equation}

If $Z_{n_0}$ has a bounded density for some $n_0$, then for all
sufficiently large $n$,
$Z_n$ have continuous bounded densities $p_n$ satisfying (\ref
{equ3.7}); see~\cite{B-RR},
Theorem~19.2. We need an extension of this theorem to the case of unbounded
densities, as well as integral variants such as (\ref{equ3.8}). The first
assertion (\ref{equ3.11})
in the next theorem is similar to the one-dimensional Theorem \ref
{Theorem3.3} in
the case
where $s = m$ is integer; cf. (\ref{equ3.5}). For the proof (which we
omit), one
may apply
Lemma~\ref{Lemma2.1} and follow the standard arguments from \cite
{B-RR}, Chapter~4.

%
\begin{Theorem}\label{Theorem3.4}
Suppose that $\E|X_1|^m < +\infty$ with some integer $m \geq2$. If, for
all sufficiently large $n$, $Z_n$ have densities $p_n$, then the
densities $\widetilde p_n$ introduced in Section~\ref{sec2} with $m_0 =
m+1$ satisfy
%
%
\begin{equation}\label{equ3.11}
\bigl(1 + |x|^m\bigr) \bigl(\widetilde p_n(x) -
\varphi_m(x) \bigr) = o \bigl(n^{-(m-2)/2} \bigr)
\end{equation}
uniformly for all $x$. In addition,
%
%
\begin{equation}\label{equ3.12}
\int\bigl(1 + |x|^m\bigr) \bigl|\widetilde p_n(x) -
\varphi_m(x)\bigr| \,dx = o \bigl(n^{-(m-2)/2} \bigr).
\end{equation}
\end{Theorem}

The second assertion is Theorem 19.5 in~\cite{B-RR}, where it is stated for
$m \geq3$
under a slightly weaker hypothesis that $X_1$ has a nonzero absolutely
continuous component. Note that, by Lemma~\ref{Lemma2.1}, it does not
matter whether
$\widetilde p_n$ or $p_n$ are used in (\ref{equ3.12}).

\section{Entropic distance to normality and moderate
deviations}\label{sec4}

Let $X_1, \break X_2,\ldots$ be independent, identically distributed random
vectors in $\R^d$
with mean zero, identity covariance matrix and such that $D(Z_n) <
+\infty$,
for all $n$ large enough.

According to Lemma~\ref{Lemma2.2} and Remark~\ref{Remark2.5}, up to an
error at most $2^{-n}$
for sufficiently large $n$,
the entropic distance to normality, $D_n = D(Z_n)$, is equal to the
relative entropy
\[
\widetilde D_n = \int\widetilde p_n(x) \log
\frac{\widetilde
p_n(x)}{\varphi(x)} \,dx,
\]
where $\varphi$ is the density of a standard normal random vector $Z$
in $\R^d$.

Given $T \geq1$, split the integral into two parts by writing
%
%
\begin{equation}\label{equ4.1}
\widetilde D_n = \int_{|x| \leq T} \widetilde
p_n(x) \log\frac{\widetilde p_n(x)}{\varphi(x)} \,dx + \int_{|x| > T}
\widetilde p_n(x) \log\frac{\widetilde p_n(x)}{\varphi
(x)} \,dx.
\end{equation}
By Theorems~\ref{Theorem3.3} and~\ref{Theorem3.4}, $\widetilde p_n$ are
uniformly bounded, that is,
$\widetilde p_n(x) \leq M$, for all $x \in\R^d$ and $n \geq1$ with some
constant $M$. Hence,\vadjust{\goodbreak} the second integral in (\ref{equ4.1}) may be
treated by virtue
of moderate deviations results (when $T$ is not too large). Indeed, since
$T \geq1$,
\[
\int_{|x| > T} \widetilde p_n(x) \log
\frac{\widetilde p_n(x)}{\varphi
(x)} \,dx \leq\int_{|x| > T} \widetilde
p_n(x) \log\frac{M}{\varphi(x)} \,dx \leq C \int_{|x| > T}
|x|^2 \widetilde p_n(x) \,dx,
\]
where $C = \frac{1}{2} + \log(1+M(2\pi)^{d/2})$.
One the other hand, using $u\log u \geq u-1$, we have a lower bound
\[
\int_{|x| > T} \widetilde p_n(x) \log
\frac{\widetilde p_n(x)}{\varphi
(x)} \,dx \geq\int_{|x| > T} \bigl(\widetilde
p_n(x) - \varphi(x) \bigr) \,dx \geq- \P\bigl\{|Z| > T\bigr\}.
\]
The two estimates give
%
%
\begin{equation}\label{equ4.2}\quad
\biggl|\int_{|x| > T} \widetilde p_n(x) \log
\frac{\widetilde
p_n(x)}{\varphi(x)} \,dx \biggr| \leq\P\bigl\{|Z| > T\bigr\} + C \int_{|x| > T}
|x|^2 \widetilde p_n(x) \,dx.
\end{equation}
This is a very general upper bound, valid for any probability
density $\widetilde p_n$ on~$\R^d$, bounded by a constant $M$ (with $C$
as above).

Following (\ref{equ4.1}), we are faced with two analytic problems. The
first one
is to give a sharp estimate of $\widetilde p_n(x) - \varphi(x)$ on a
relatively large Euclidean ball \mbox{$|x| \leq T$}. Clearly, $T$ has to be
small enough, so that results like local limit theorems, such as
Theorems~\ref{Theorem3.2}--\ref{Theorem3.4} may be applied. The second
problem is to give a sharp upper bound of the last integral in (\ref{equ4.2}).
To this aim, we need moderate deviations inequalities, so that Theorems
\ref{Theorem3.1} and~\ref{Theorem3.4} are applicable. Anyway, in order
to use both types of results we are forced to choose $T$ from a very
narrow window only. This value turns out to be approximately
%
%
\begin{equation}\label{equ4.3}
T_n = \sqrt{(s-2)\log n + s \log\log n + \rho_n}\qquad (s>2),
\end{equation}
where $\rho_n \rightarrow+\infty$ is a sufficiently slowly growing sequence
(whose growth will be restricted by the decay of the $n$-dependent
constants in
$o$-expressions of Theorems~\ref{Theorem3.2}--\ref{Theorem3.4}). In the
case $s=2$, one may put
$T_n = \sqrt{\rho_n}$
such that $T_n \rightarrow+\infty$ is a sufficiently slowly growing sequence.

%
\begin{Lemma}[(The case $d=1$ and $s$ real)]\label{Lemma4.1}
If $\E X_1 = 0$,
$\E X_1^2 = 1$,
$\E|X_1|^s < + \infty$ $(s \geq2)$, then
%
%
\begin{equation}\label{equ4.4}
\int_{|x| > T_n} x^2 \widetilde p_n(x) \,dx
= o \bigl((n \log n)^{-(s-2)/2} \bigr).
\end{equation}
\end{Lemma}

%
\begin{Lemma}[(The case $d \geq2$ and $s$ integer)]\label{Lemma4.2}
If $X_1$ has mean zero and identity covariance matrix, and
$\E|X_1|^m < + \infty$, then
%
%
\begin{equation}\label{equ4.5}
\int_{|x| > T_n} x^2 \widetilde p_n(x) \,dx
= o \bigl(n^{-(m-2)/2} (\log n)^{-(m-d)/2} \bigr)\qquad (m \geq3)
\end{equation}
and $\int_{|x| > T_n} x^2 \widetilde p_n(x) \,dx = o(1)$ in the case $m=2$.
\end{Lemma}

Note that plenty of results and techniques concerning moderate deviations
have been developed by now. Useful estimates can be found, for example,
in~\cite{G-H}.
Restricting ourselves to integer values of $s = m$, one may argue as follows.

\begin{pf*}{Proof of Lemma~\ref{Lemma4.2}}
Given $T \geq1$, write
%
%
\begin{eqnarray}\label{equ4.6}
\int_{|x| > T} |x|^2 \widetilde
p_n(x) \,dx & \leq& \frac{1}{T^{m-2}} \int|x|^m
\widetilde p_n(x) \,dx
\nonumber\\
& \leq & \frac{1}{T^{m-2}} \int|x|^m \bigl|\widetilde
p_n(x) - \varphi_m(x)\bigr| \,dx \\
&&{}+ \frac{1}{T^{m-2}} \int
_{|x| > T} |x|^m \varphi_m(x) \,dx.\nonumber
\end{eqnarray}
By Theorem~\ref{Theorem3.4} [cf. (\ref{equ3.12})] the first integral in
(\ref{equ4.6}) is bounded by
$o(n^{-(m-2)/2})$.

From the definition of $q_k$ it follows that $q_k(x) = N(x)\varphi(x)$
with some polynomial $N$ of degree at most $3(m-2)$; cf. Section \ref
{sec6} for details.
Hence, from (\ref{equ3.1}), $\varphi_m(x) \leq2 \varphi(x)$ on the
balls of
large radii
$|x| < n^\delta$ with sufficiently large $n$ (where $0 < \delta<
\frac{1}{2}$).
On the other hand, with some constants $C_d, C_d'$ depending on the
dimension only,
%
%
\begin{equation}\label{equ4.7}\qquad
\int_{|x| > T} |x|^m \varphi(x) \,dx =
C_d \int_T^{+\infty} r^{m + d - 1}
e^{-r^2/2} \,dr \leq C_d' T^{m+d-2}
e^{-T^2/2}.
\end{equation}
But for $T=T_n$ and $s = m \geq3$, we have
$e^{-T^2/2} = T^{-m} o(n^{-(m-2)/2})$, so by (\ref{equ4.6}) and (\ref{equ4.7}),
\[
\int_{|x| > T_n} |x|^2 \widetilde p_n(x)
\,dx \leq C \biggl(\frac{1}{T^{m-2}} + \frac{1}{T^{m-d}} \biggr) o
\bigl(n^{-(m-2)/2}\bigr).
\]
Since $T_n$ is of order $\sqrt{\log n}$, (\ref{equ4.5}) follows.
Furthermore, in the case $m=2$, (\ref{equ4.6}) gives the desired relation
\[
\int_{|x| > T_n} |x|^2 \widetilde p_n(x)
\,dx \leq o(1) + \int_{|x| > T_n} |x|^2 \varphi(x) \,dx
\rightarrow0 \qquad(n \rightarrow\infty).
\]
\upqed\end{pf*}

\begin{pf*}{Proof of Lemma~\ref{Lemma4.1}}
The above argument also works for $d=1$, but it can be refined applying
Theorem~\ref{Theorem3.1} for real $s$. The case $s=2$ is already
covered, so let
$s>2$.

In view of decomposition (\ref{equ2.5}), integrating by parts, we have,
for any
$T \geq0$,
%
%
\begin{eqnarray}\quad
\label{equ4.8}
&&
(1 - \ep_n)\int_{|x| > T} x^2
\widetilde p_n(x) \,dx \nonumber\\[-8pt]\\[-8pt]
&&\qquad \leq \int_{|x| > T}
x^2 p_n(x) \,dx
= \int_{|x| > T}
x^2 \,dF_n(x)\nonumber
\\
\label{equ4.9}
&&\qquad= T^2 \bigl(1-F_n(T) +
F_n(-T) \bigr)
+ 2 \int_T^{+\infty} x
\bigl(1-F_n(x) + F_n(-x) \bigr) \,dx,
\end{eqnarray}
where $F_n$ denotes the distribution function of $Z_n$. [Note that the first
inequality in (\ref{equ4.8}) should be just ignored in the case, where
$p$ is bounded.]

By (\ref{equ3.3}),
\[
F_n(x) = \Phi_m(x) + \frac{r_n(x)}{n^{(s-2)/2}}
\frac{1}{1 + |x|^s},\qquad r_n = \sup_x \bigl|r_n(x)\bigr|
\rightarrow0 \qquad(n \rightarrow\infty).
\]
Hence, the first term in (\ref{equ4.9}) can be replaced with
%
%
\begin{equation}\label{equ4.10}
T^2 \bigl(1-\Phi_m(T) + \Phi_m(-T) \bigr)
\end{equation}
at the expense of an error not exceeding (for the values $T \sim\sqrt
{\log n}$)
%
%
\begin{equation}\label{equ4.11}
\frac{2r_n}{n^{(s-2)/2}} \frac{T^2}{1 + T^s} = o \bigl((n \log n)^{-(s-2)/2}
\bigr).
\end{equation}
Similarly, the integral in (\ref{equ4.9}) can be replaced with
%
%
\begin{equation}\label{equ4.12}
\int_T^{+\infty} x \bigl(1-\Phi_m(x)
+ \Phi_m(-x) \bigr) \,dx
\end{equation}
at the expense of an error not exceeding
%
%
\begin{equation}\label{equ4.13}
\frac{2r_n}{n^{(s-2)/2}} \int_T^{+\infty}
\frac{x \,dx}{1 + x^s} = o \bigl((n \log n)^{-(s-2)/2} \bigr).
\end{equation}

To explore the behavior of expressions (\ref{equ4.10}) and
(\ref{equ4.12}) for $T = T_n$ using precise asymptotics as in (\ref{equ4.3}),
recall that, by (\ref{equ3.2}),
\[
1 - \Phi_m(x) = 1 - \Phi(x) - \sum_{k=1}^{m-2}
Q_k(x) n^{-k/2}.
\]
Moreover, we note that $Q_k(x) = N_{3k-1}(x) \varphi(x)$, where $N_{3k-1}$
is a polynomial of degree\vspace*{1pt} at most $3k-1$. Thus these functions admit a bound
$|Q_k(x)| \leq C_m (1 + |x|^{3m}) \varphi(x)$ with some constants $C_m$
(depending on $m$ and the cumulants $\gamma_3,\ldots,\gamma_m$ of $X_1$),
which implies with some other constants
%
%
\begin{equation}\label{equ4.14}
\bigl|1 - \Phi_m(x)\bigr| \leq\bigl(1 - \Phi(x)\bigr) + \frac{C_m (1 +
|x|^{3m})}{\sqrt{n}}
\varphi(x).
\end{equation}
Hence, using $1 - \Phi(x) < \frac{\varphi(x)}{x}$ ($x > 0$), we get
%
%
\begin{eqnarray}\label{equ4.15}
T_n^2 \bigl|1 - \Phi_m(T_n)\bigr| &\leq& C
T_n^2 \bigl(1 - \Phi(T_n)\bigr) \leq C
T_n e^{-T_n^2/2} \nonumber\\[-8pt]\\[-8pt]
&=& o \bigl((n \log n)^{-(s-2)/2} \bigr).\nonumber
\end{eqnarray}
A similar bound also holds for $T_n^2 |\Phi_m(-T_n)|$.\vadjust{\goodbreak}

Now, we use (\ref{equ4.14}) to estimate (\ref{equ4.12}) with $T = T_n$
up to a constant by
\[
\int_T^\infty x \bigl(1 - \Phi(x)\bigr) \,dx < 1
- \Phi(T) = o \bigl((n \log n)^{-(s-2)/2} \bigr).
\]

It remains to combine the last relation with (\ref{equ4.11}), (\ref
{equ4.13}) and (\ref{equ4.15}).
Since $\ep_n \rightarrow0$ in (\ref{equ4.8}), Lemma~\ref{Lemma4.1} follows.
\end{pf*}

%
\begin{Remark}\label{Remark4.3}
Note that the probabilities $\P\{|Z| > T\}$ appearing
in (\ref{equ4.2}) yield a smaller contribution for $T = T_n$ in
comparison with
the right-hand
sides of (\ref{equ4.4}) and (\ref{equ4.5}). Indeed, we have
$\P\{|Z| > T\} \leq C_d T^{d-2} e^{-T^2/2}$ (\mbox{$T \geq1$}).
Hence, relations (\ref{equ4.4}) and (\ref{equ4.5}) may be extended to
the integrals
\[
\int_{|x| > T_n} \widetilde p_n(x) \log
\frac{\widetilde p_n(x)}{\varphi
(x)} \,dx.
\]
\end{Remark}

\section{Taylor-type expansion for the entropic distance}\label{sec5}

In this section we provide the last auxiliary step toward the proof of
Theorem~\ref{Theorem1.1}.
In order to describe the multidimensional case, let $X_1,X_2,\ldots$ be
independent
identically distributed random vectors in $\R^d$ with mean zero,
identity covariance matrix, and such that $D(Z_{n_0}) < +\infty$ for
some $n_0$.

If $p_{n_0}$ is bounded, then the densities $p_n$ of $Z_n$ $(n \geq n_0)$
are uniformly bounded, and we put $\widetilde p_n = p_n$. Otherwise, we
use the
modified densities $\widetilde p_n$ according to the construction of
Section~\ref{sec2}.
In particular, if $\widetilde Z_n$ has density~$\widetilde p_n$, then
$|D(\widetilde Z_n\|Z) - D(Z_n)| < 2^{-n}$ for all $n$ large enough
(where $Z$ is a standard normal random vector; cf. Lemma~\ref{Lemma2.2}
and Remark~\ref{Remark2.5}).
Moreover, by Lemmas~\ref{Lemma4.1},~\ref{Lemma4.2} and Remark~\ref{Remark4.3},
%
%
\begin{equation}\label{equ5.1}
\biggl|D(Z_n) - \int_{|x| \leq T_n} \widetilde
p_n(x) \log\frac{\widetilde p_n(x)}{\varphi(x)} \,dx \biggr| = o(\Delta_n),
\end{equation}
where $T_n$ are defined in (\ref{equ4.3}) and
%
%
\begin{equation}\label{equ5.2}
\Delta_n = n^{-(s-2)/2} (\log n)^{-(s-\max(d,2))/2}
\end{equation}
(with the convention that $\Delta_n = 1$ for the critical case $s=2$).

Thus, all information about the asymptotics of $D(Z_n)$ is contained in
the integral in (\ref{equ5.1}). More precisely, writing a Taylor
expansion for
$\widetilde p_n$ using the approximating functions $\varphi_m$ in
Theorems~\ref{Theorem3.2}--\ref{Theorem3.4} leads to the following
representation (which is more convenient in applications such as
Corollary~\ref{Corollary1.2}).

%
\begin{Theorem}\label{Theorem5.1} Let $\E|X_1|^s < +\infty$ $(s \geq2)$,
assuming that
$s$ is integer in case $d \geq2$. Then
%
%
\begin{eqnarray}\label{equ5.3}
D(Z_n) &=& \sum_{k=2}^{m-2}
\frac{(-1)^k}{k(k-1)} \int\bigl(\varphi_m(x) - \varphi(x)
\bigr)^k \,\frac{dx}{\varphi(x)^{k-1}} \nonumber\\[-8pt]\\[-8pt]
&&{}+ o(\Delta_n) \qquad\bigl(m = [s]
\bigr).\nonumber
\end{eqnarray}
\end{Theorem}

Note that in the case $2 \leq s < 4$ there are no expansion terms in
the sum of (\ref{equ5.3}) which
then simplifies to $D(Z_n) = o(\Delta_n)$.

\begin{pf*}{Proof of Theorem~\ref{Theorem5.1}}
In terms of $L(u) = u\log u$, rewrite the integral in
(\ref{equ5.1}) as
%
%
\begin{eqnarray}\label{equ5.4}
\widetilde D_{n,1} &=& \int_{|x| \leq T_n} L \biggl(
\frac{\widetilde
p_n(x)}{\varphi(x)} \biggr) \varphi(x) \,dx \nonumber\\[-8pt]\\[-8pt]
&=& \int_{|x| \leq T_n} L
\bigl(1 + u_m(x) + v_n(x) \bigr) \varphi(x) \,dx,\nonumber
\end{eqnarray}
where
\[
u_m(x) = \frac{\varphi_m(x) - \varphi(x)}{\varphi(x)},\qquad v_n(x) =
\frac{\widetilde p_n(x) - \varphi_m(x)}{\varphi(x)}.
\]

By Theorems~\ref{Theorem3.3} and~\ref{Theorem3.4}, more precisely, by
(\ref{equ3.6}) for $d=1$, and by (\ref{equ3.11}) for $d \geq2$ and $s=m$
integer, in the region $|x| = O(n^\delta)$ with an appropriate
$\delta>0$, we have
%
%
\begin{equation}\label{equ5.5}
\bigl|\widetilde p_n(x) - \varphi_m(x)\bigr| \leq
\frac{r_n}{n^{(s-2)/2}} \frac
{1}{1 + |x|^s},\qquad r_n \rightarrow0.
\end{equation}
Since $\varphi(x) (1 + |x|^s)$ is decreasing as a function of $|x|$
for large $|x|$, we obtain, for all $|x| \leq T_n$,
\[
\bigl|v_n(x)\bigr| \leq C \frac{r_n}{n^{(s-2)/2}} \frac{e^{T_n^2/2}}{T_n^s}
\leq
C' r_n e^{\rho_n/2}.
\]
The last expression tends to zero by a suitable choice of $\rho_n
\rightarrow\infty$
which we will assume from now on.
In particular, for $n$ large enough, $|v_n(x)| < \frac{1}{4}$ in $|x|
\leq T_n$.

From the definitions of $q_k$ and $\varphi_m$ [cf. (\ref{equ1.2}), (\ref
{equ3.1}) and (\ref{equ3.10})],
it follows that
%
%
\begin{equation}\label{equ5.6}
\bigl|u_m(x)\bigr| \leq C_m \frac{1 + |x|^{3(m-2)}}{\sqrt{n}}
\end{equation}
with some constants depending on $m$ and the cumulants, only. Thus, we
also have
$|u_m(x)| < \frac{1}{4}$ for $|x| \leq T_n$ with sufficiently large $n$.

Now, by Taylor's formula, for $|u| \leq\frac{1}{4}$, $|v| \leq\frac{1}{4}$,
\[
L(1 + u + v) = L(1+u) + v + 2\theta_1 u v + \theta_2
v^2
\]
with some $|\theta_j| \leq1$ depending on $(u,v)$. Applying this approximation
with $u = u_m(x)$ and $v = v_n(x)$, we see that $v_n(x)$ can be removed
from the right-hand side of (\ref{equ5.4}) at the expense of an error
not exceeding
$|J_1| + J_2 + J_3$, where
\[
J_1 = \int_{|x| \leq T_n} \bigl(\widetilde
p_n(x) - \varphi_m(x) \bigr) \,dx,\qquad J_2 = \int
_{|x| \leq T_n} \bigl|u_m(x)\bigr| \bigl|\widetilde p_n(x) -
\varphi_m(x)\bigr| \,dx
\]
and
\[
J_3 = \int_{|x| \leq T_n} \frac{(\widetilde p_n(x) - \varphi
_m(x))^2}{\varphi(x)} \,dx.
\]

But
%
%
\begin{eqnarray}\label{equ5.7}
|J_1| &=& \biggl|\int_{|x| > T_n} \bigl(\widetilde
p_n(x) - \varphi_m(x)\bigr) \,dx \biggr| \nonumber\\[-8pt]\\[-8pt]
&\leq&\int
_{|x| > T_n} \widetilde p_n(x) \,dx + \int
_{|x| > T_n} \bigl|\varphi_m(x)\bigr| \,dx.\nonumber
\end{eqnarray}
By Lemmas~\ref{Lemma4.1} and~\ref{Lemma4.2}, the first integral on the
right-hand side is $T_n^2$-times
smaller than $o(\Delta_n)$. Also, by (\ref{equ5.6}), the last integral
in (\ref{equ5.7})
is bounded by
\begin{eqnarray*}
&&
\int_{|x| > T_n} \bigl|\varphi_m(x) -
\varphi(x)\bigr| \,dx + \int_{|x| > T_n} \varphi(x) \,dx
\\
&&\qquad \leq\frac{C_m}{\sqrt{n}} \int_{|x| > T_n}
\bigl(1 + |x|^{3(m-2)}\bigr) \varphi(x) \,dx + \P\bigl\{|Z| > T_n\bigr\}
= o(\Delta_n).
\end{eqnarray*}
As a result, $J_1 = o(\Delta_n)$.

Applying (\ref{equ5.6}) once more and then relation (\ref{equ3.12}), we
may also conclude that
\[
J_2 \leq C_m \frac{1 + T_n^{3(m-2)}}{\sqrt{n}} \int
_{|x| \leq T_n} \bigl|\widetilde p_n(x) - \varphi_m(x)\bigr|
\,dx = o(\Delta_n).
\]

Finally, using (\ref{equ5.5}) with $s>2$, we get, up to some constants,
\begin{eqnarray*}
J_3 & \leq& C \frac{r_n^2}{n^{s-2}} \int_{|x| \leq T_n}
\frac{e^{|x|^2/2}}{1 +
|x|^{2s}} \,dx \leq C_d \frac{r_n^2}{n^{s-2}} \int
_1^{T_n} r^{d-2s-1} e^{r^2/2} \,dr
\\
& \leq& C_d' \frac{r_n^2}{n^{s-2}} \frac{1}{T_n^{2s-d+2}}
e^{T_n^2/2} = o \biggl(\frac{1}{n^{(s-2)/2} (\log n)^{(s-d+2)/2}}
\biggr) = o(\Delta_n).
\end{eqnarray*}
If $s=2$, all these steps are valid as well and give
\[
J_3 \leq C_d' \frac{r_n^2}{n^{s-2}}
\frac{1}{T_n^{2s-d+2}} e^{T_n^2/2} \rightarrow0
\]
for a suitably chosen $T_n \rightarrow+\infty$.

Thus, at the expense of an error not exceeding $o(\Delta_n)$ one may
remove $v_n(x)$ from (\ref{equ5.4}), and we obtain the relation
%
%
\begin{equation}\label{equ5.8}
\widetilde D_{n,1} = \int_{|x| \leq T_n} L \bigl(1 +
u_m(x) \bigr) \varphi(x) \,dx + o(\Delta_n),
\end{equation}
which contains specified expansion terms, only.

Moreover, $u_m(x) = u_2(x) = 0$ for $2 \leq s < 3$, and then the
theorem is~proved.\looseness=-1

Next, we consider the case $s \geq3$. By Taylor's expansion around
zero, we get, whenever
$|u| < \frac{1}{4}$, for some
positive constants $\theta_m$,
\[
L(1 + u) = u + \sum_{k=2}^{m-2}
\frac{(-1)^k}{k(k-1)} u^k +\theta u^{m-1},\qquad |\theta| \leq
\theta_m,
\]
assuming that the sum has no terms in the case $m=3$. Hence, with some
$|\theta| \leq \theta_m$,
%
%
\begin{eqnarray}
\label{equ5.9}
&&\int_{|x| \leq T_n} L \bigl(1 + u_m(x) \bigr) \varphi(x)
\,dx \nonumber\\[-9pt]\\[-9pt]
&&\qquad = \int_{|x| \leq T_n} \bigl(\varphi_m(x) -
\varphi(x) \bigr) \,dx\nonumber
\\[-2pt]
\label{equ5.10}
&&\qquad\quad{} + \sum_{k=2}^{m-2}
\frac{(-1)^k}{k(k-1)} \int_{|x| \leq T_n} u_m(x)^k
\varphi(x) \,dx\nonumber\\[-9pt]\\[-9pt]
&&\qquad\quad{} +\theta\int_{\R^d} \bigl|u_m(x)\bigr|^{m-1}
\varphi(x) \,dx.\nonumber
\end{eqnarray}

For $n$ large enough, by (\ref{equ5.6}), the second integral in (\ref
{equ5.9}) has an
absolute value
\[
\biggl|\int_{|x| > T_n} \bigl(\varphi_m(x) - \varphi(x)
\bigr) \,dx \biggr| \leq\frac{C}{\sqrt{n}} \int_{|x| > T_n} \bigl(1 +
|x|^{3(m-2)}\bigr) \varphi(x) \,dx = o(\Delta_n).
\]
This proves the theorem in the case $3 \leq s < 4$ (when $m=3$).

Now, let $s \geq4$. The last integral in (\ref{equ5.10}) can be
estimated again
by virtue of (\ref{equ5.6}) by
\[
\frac{C}{n^{(m-1)/2}} \int_{\R^d} \bigl(1 + |x|^{3(m-1)(m-2)}
\bigr) \varphi(x) \,dx = o(\Delta_n).
\]
In addition, the first integral in (\ref{equ5.10}) can be extended to
the whole space
at the expense of an error not exceeding (for all $n$ large enough)
\begin{eqnarray*}
\int_{|x| > T_n} \bigl|u_m(x)\bigr|^k \varphi(x) \,dx
& \leq& \frac{C}{n^{k/2}} \int_{|x| > T_n} \bigl(1 +
|x|^{3k(m-2)} \bigr) \varphi(x) \,dx
\\
& \leq& \frac{C' T_n^{3k(m-2)}}{\sqrt{n}} e^{-T_n^2/2} = o(\Delta_n).
\end{eqnarray*}

Collecting these estimates in (\ref{equ5.9}) and (\ref{equ5.10}) and
applying them in (\ref{equ5.8}),
we arrive at
\[
\widetilde D_{n,1} = \sum_{k=2}^{m-2}
\frac{(-1)^k}{k(k-1)} \int u_m(x)^k \varphi(x) \,dx + o(
\Delta_n).
\]
It remains to apply (\ref{equ5.1}). Thus, Theorem~\ref{Theorem5.1} is proved.\vadjust{\goodbreak}
\end{pf*}

\section{\texorpdfstring{Theorem \protect\ref{Theorem1.1} and its multidimensional extension}
{Theorem 1.1 and its multidimensional extension}}\label{sec6}

The desired representation (\ref{equ1.3}) of Theorem~\ref{Theorem1.1}
can be deduced from
Theorem~\ref{Theorem5.1}. Note that the latter covers the
multidimensional case as well,
although under somewhat stronger moment assumptions.

Thus, let $(X_n)_{n \geq1}$ be independent identically distributed random
vectors in $\R^d$ with finite second moment. If the normalized sum
$Z_n = (X_1 + \cdots+ X_n)/\sqrt{n}$ has density $p_n(x)$, the entropic
distance to Gaussianity is defined as in dimension one to be the
relative entropy
\[
D(Z_n) = \int p_n(x) \log\frac{p(x)}{\varphi_{a,\Sigma}(x)} \,dx
\]
with respect to the normal law on $\R^d$ with the same mean $a = \E
X_1$ and
covariance matrix $\Sigma= \Var(X_1)$. This quantity is affine invariant,
and in this sense it does not depend on $(a,\Sigma)$.

%
\begin{Theorem}\label{Theorem6.1}
If $D(Z_{n_0}) < +\infty$ for some $n_0$, then $D(Z_n) \rightarrow0$,
as \mbox{$n \rightarrow\infty$}. Moreover, given that $\E|X_1|^s < +\infty$
$(s \geq2)$, and that $X_1$ has mean zero and identity covariance
matrix, we have
%
%
\begin{equation}\label{equ6.1}
D(Z_n)= \frac{c_1}{n} + \frac{c_2}{n^2} + \cdots+
\frac{c_{[(m-2)/2]}}{n^{[(m-2)/2]}} + o (\Delta_n) \qquad\bigl(m =
[s]\bigr),
\end{equation}
where $\Delta_n$ are defined in (\ref{equ5.2}), and where we assume
that $s$
is integer
in case $d \geq2$.
\end{Theorem}

Here, as in Theorem~\ref{Theorem1.1}, each coefficient $c_j$ is defined
according to
(\ref{equ1.4}) again.
It may be represented as a certain polynomial in the cumulants
$\gamma_\nu$, $3 \leq|\nu| \leq 2j+1$.

\begin{pf*}{Proof of Theorem~\ref{Theorem6.1}}
We shall start from the representation (\ref{equ5.3}) of Theorem \ref
{Theorem5.1}, so let us
return to definition (\ref{equ3.1}),
\[
\varphi_m(x) - \varphi(x) = \sum_{r=1}^{m-2}
q_r(x) n^{-r/2}.
\]
In the case $2 \leq s < 3$ (i.e., for $m=2$), the right-hand side
contains no terms
and is therefore vanishing. Anyhow, raising this sum to the power $k
\geq2$ leads to
\[
\bigl(\varphi_m(x) - \varphi(x) \bigr)^k = \sum
_j n^{-j/2} \sum
q_{r_1}(x) \cdots q_{r_k}(x),
\]
where the inner sum is carried out over all positive integers
$r_1,\ldots,r_k \leq m-2$ such that $r_1 + \cdots+ r_k = j$.
Respectively, the $k$th integral in (\ref{equ5.3}) is equal to
%
%
\begin{equation}\label{equ6.2}
\sum_{j} n^{-j/2} \sum\int
q_{r_1}(x) \cdots q_{r_k}(x) \frac
{dx}{\varphi(x)^{k-1}}.
\end{equation}

Here the integrals are vanishing for odd $j$. In dimension one, this follows
directly from definition (\ref{equ1.2}) of $q_r$ and the following
property of
the Chebyshev--Hermite polynomials~\cite{Sz}
%
%
\begin{equation}\label{equ6.3}
\int_{-\infty}^{+\infty} H_{r_1}(x) \cdots
H_{r_k}(x) \varphi(x) \,dx = 0 \qquad(r_1 + \cdots+
r_k \mbox{ is odd}).
\end{equation}
As for the general case, let us look at the structure of the functions
$q_r$. Given a multi-index $\nu= (\nu_1,\ldots,\nu_d)$ with integers
$\nu_1,\ldots,\nu_d \geq1$, define
$H_\nu(x_1,\ldots,\break x_d) = H_{\nu_1}(x_1) \cdots H_{\nu_d}(x_d)$, so that
\[
\int e^{i\langle t,x\rangle} H_\nu(x) \varphi(x) \,dx = (it)^\nu
e^{-|t|^2/2},\qquad t \in\R^d.
\]
Hence, by definition (\ref{equ3.10}),
%
%
\begin{equation}\label{equ6.4}
q_r(x) = \varphi(x) \sum_\nu
a_\nu H_\nu(x),
\end{equation}
where the coefficients $a_\nu$ emerge from the expansion $P_r(it) =
\sum_\nu a_\nu(it)^\nu$.
Using (\ref{equ3.9}), write these polynomials as
%
%
\begin{equation}\label{equ6.5}
P_r(it) = \sum\frac{1}{l_1!\cdots l_r!} \biggl(\sum
_{|\nu| = 3} \gamma_\nu\frac{(it)^\nu}{\nu!}
\biggr)^{l_1} \cdots\biggl(\sum_{|\nu| = r+2}
\gamma_\nu\frac{(it)^\nu}{\nu!} \biggr)^{l_r},
\end{equation}
where the outer summation is performed over all nonnegative integer solutions
$(l_1,\ldots,l_r)$ to the equation $l_1 + 2 l_2 + \cdots+ r l_r = r$.
Removing the brackets of the inner sums, we obtain a linear combination
of the power polynomials $(it)^\nu$ with exponents of order
%
%
\begin{equation}\label{equ6.6}
|\nu| = 3l_1 + \cdots+ (r+2) l_r = r + 2
b_l,\qquad b_l = l_1 + \cdots+
l_r.
\end{equation}
In particular, $r+2 \leq|\nu| \leq3r$, so that $P_r(it)$ is a
polynomial of degree
at most~$3r$, and thus $\varphi_m(x) = N(x) \varphi(x)$, where
$N(x)$ is a polynomial of degree at most \mbox{$3(m-2)$}.

Moreover, from (\ref{equ6.4}) and (\ref{equ6.6}) it follows that
%
%
\begin{equation}\label{equ6.7}
\frac{q_{r_1}(x) \cdots q_{r_k}(x)}{\varphi(x)^{k-1}} = \varphi(x)
\sum a_{\nu^{(1)}} \cdots
a_{\nu^{(k)}} H_{\nu^{(1)}}(x) \cdots H_{\nu^{(k)}}(x),
\end{equation}
where $|\nu^{(1)}| + \cdots+ |\nu^{(k)}| = r_1 + \cdots+ r_k
(\operatorname{mod} 2)$.
Hence, if $r_1 + \cdots+ r_k$ is odd, the sum
\[
\bigl|\nu^{(1)}\bigr| + \cdots+ \bigl|\nu^{(k)}\bigr| = \sum
_{i=1}^d \bigl(\bigl|\nu_i^{(1)}\bigr|
+ \cdots+ \bigl|\nu_i^{(k)}\bigr| \bigr)
\]
is odd as well. But then at least one of the inner sums, say with
coordinate~$i$, must be odd as well. Hence in this case, the integral
of (\ref{equ6.7})
over $x_i$ will vanish by property (\ref{equ6.3}).

Thus, in expression (\ref{equ6.2}), only even values of $j$ should be
taken into
account.\vadjust{\goodbreak}

Moreover, since the terms containing $n^{-j/2}$ with $j > s-2$ will be
absorbed into the remainder $\Delta_n$ in relation (\ref{equ6.1}), we
get from (\ref{equ5.3}) and (\ref{equ6.2}),
\[
D(Z_n) = \sum_{k=2}^{m-2}
\frac{(-1)^k}{k(k-1)} \sum_{\mathrm{even} j=2}^{m-2}
n^{-j/2} \sum\int q_{r_1}(x) \cdots
q_{r_k}(x) \,\frac{dx}{\varphi(x)^{k-1}} + o(\Delta_n).
\]
Replace now $j$ with $2j$ and rearrange the summation. Then
\[
D(Z_n) = \sum_{2j \leq m-2} \frac{c_j}{n^{j}} + o(\Delta_n)
\]
with
\[
c_j = \sum_{k=2}^{m-2}
\frac{(-1)^k}{k(k-1)} \sum\int q_{r_1}(x) \cdots
q_{r_k}(x) \,\frac{dx}{\varphi(x)^{k-1}}.
\]
Here the inner summation is carried out over all positive integers
$r_1,\ldots,r_k \leq m-2$ such that $r_1 + \cdots+ r_k = 2j$. This implies
$k \leq2j$. Furthermore, $2j \leq m-2$ is equivalent to $j \leq[\frac
{s-2}{2}]$.
As a result, we arrive at the required relation (\ref{equ6.1}) with
%
%
\begin{equation}\label{equ6.8}
c_j = \sum_{k=2}^{2j}
\frac{(-1)^k}{k(k-1)} \sum_{r_1 + \cdots+ r_k
= 2j} \int
q_{r_1}(x) \cdots q_{r_k}(x) \,\frac{dx}{\varphi(x)^{k-1}}.
\end{equation}
Thus, Theorem~\ref{Theorem6.1} and therefore Theorem~\ref{Theorem1.1}
are proved.
\end{pf*}

\begin{Remark*} In order to show that $c_j$ is a polynomial in the cumulants
$\gamma_\nu$, $3 \leq|\nu| \leq2j+1$, first note that
$r_1 + \cdots+ r_k = 2j$, $r_1,\ldots,r_k \geq1$ imply
$2j \geq\max_i r_i + (k-1)$, so $\max_i r_i \leq2j-1$.
Thus, the maximal index for the functions $q_{r_i}$ in (\ref{equ6.8})
does not
exceed $2j-1$. On the other hand, it follows from (\ref{equ6.4}) and
(\ref{equ6.5}) that
$P_r$ and $q_r$ are polynomials in the same set of the cumulants; more
precisely,
$P_r$ is a polynomial in $\gamma_\nu$ with $3 \leq|\nu| \leq r+2$.
\end{Remark*}

\begin{pf*}{Proof of Corollary~\ref{Corollary1.2}}
By Theorem~\ref{Theorem5.1} [cf. (\ref{equ5.3})],
%
%
\begin{equation}\label{equ6.9}
D(Z_n) = \sum_{k=2}^{m-2}
\frac{(-1)^k}{k(k-1)} \int\bigl(\varphi_m(x) - \varphi(x)
\bigr)^k \,\frac{dx}{\varphi(x)^{k-1}} + o (\Delta_n ).
\end{equation}
Assume that $m \geq4$ and $\gamma_3 = \cdots= \gamma_{k-1} = 0$ for
a given integer $3 \leq k \leq m$. (This is no restriction, when $k = 3$.)
Then, by (\ref{equ1.2}), $q_1 = \cdots= q_{k-3} = 0$, while
$q_{k-2}(x) = \frac{\gamma_k}{k!} H_k(x) \varphi(x)$. Hence,
according to
definition (\ref{equ3.1}),
\[
\varphi_m(x) - \varphi(x) = \frac{\gamma_k}{k!} H_k(x)
\varphi(x) \frac{1}{n^{(k-2)/2}} + \sum_{j=k-1}^{m-2}
\frac{q_j(x)}{n^{j/2}},
\]
where the sum is empty in the case $m=3$.
Therefore, the sum in (\ref{equ1.3}) will contain powers of $1/n$
starting from
$1/n^{k-2}$, and the leading coefficient is due to the quadratic term
in (\ref{equ6.9})\vadjust{\goodbreak}
when $k=2$. More precisely, if $k-2 \leq\frac{m-2}{2}$, we get that
$c_1 = \cdots= c_{k-3} = 0$, and
%
%
\begin{equation}\label{equ6.10}
c_{k-2} = \frac{\gamma_k^2}{2 k!^{ 2}} \int_{-\infty}^{+\infty}
H_k(x)^2 \varphi(x) \,dx = \frac{\gamma_k^2}{2 k!}.
\end{equation}
Hence, if $k \leq\frac{m}{2}$, (\ref{equ6.9}) yields
$D(Z_n) = \frac{\gamma_k^2}{2 k!} \frac{1}{n^{k-2}} + O
(n^{-(k-1)} )$.
Otherwise, the $O$-term should be replaced by $o((n \log n)^{-(s-2)/2})$.
Thus Corollary~\ref{Corollary1.2} is proved.
\end{pf*}

By a similar argument, the conclusion may be extended to the multidimensional
case. Indeed, if $\gamma_\nu= 0$, for all $3 \leq|\nu| < k$, then
by (\ref{equ6.5}),
$P_1 = \cdots= P_{k-3} = 0$, while
\[
P_{k-2}(it) = \sum_{|\nu| = k}
\gamma_\nu\frac{(it)^\nu}{\nu!}.
\]
Correspondingly, in (\ref{equ6.4}) we have $q_1 = \cdots= q_{k-3} = 0$ and
$q_{k-2}(x) = \varphi(x)\* \sum_{|\nu| = k} \frac{\gamma_\nu}{\nu
!} H_\nu(x)$.
Therefore,
\[
\varphi_m(x) - \varphi(x) = \varphi(x) \sum
_{|\nu| = k} \frac{\gamma_\nu}{\nu!} H_\nu(x)
\frac
{1}{n^{(k-2)/2}} + \sum_{j=k-1}^{m-2}
\frac{q_j(x)}{n^{j/2}}.
\]
Applying this relation in (\ref{equ6.9}), we arrive at (\ref{equ6.1}) with
$c_1 = \cdots= c_{k-3} = 0$
and, by orthogonality of the polynomials $H_\nu$,
\[
c_{k-2} = \frac{1}{2} \int\biggl(\sum
_{|\nu| = k} \frac{\gamma_\nu}{\nu!} H_\nu(x)
\biggr)^2 \varphi(x) \,dx = \frac{1}{2} \sum
_{|\nu| = k} \frac{\gamma_\nu^2}{\nu!}.
\]
We may summarize our findings as follows.

%
\begin{Corollary}\label{Corollary6.2}
Let $(X_n)_{n \geq1}$ be i.i.d. random vectors in $\R^d$ $(d \geq2)$
with mean zero and identity covariance matrix. Suppose that $\E|X_1|^m
< +\infty$, for some integer $m \geq4$, and $D(Z_{n_0}) < +\infty$, for
some $n_0$. Given $k = 3,4,\ldots,m$, if $\gamma_\nu= 0$ for all $3
\leq|\nu| < k$, we have
%
%
\begin{equation}\label{equ6.11}\qquad
D(Z_n) = \frac{1}{2 n^{k-2}} \sum_{|\nu| = k}
\frac{\gamma_\nu^2}{\nu!} + O \biggl(\frac{1}{n^{k-1}} \biggr) +
o \biggl(
\frac{1}{n^{(m-2)/2} (\log n)^{(m-d)/2} } \biggr).
\end{equation}
\end{Corollary}

The conclusion corresponds to Corollary~\ref{Corollary1.2}, if we
replace $d$ with $2$
in the remainder on the right-hand side.

As in dimension one, when $\E X_1^{2k} < +\infty$, the $o$-term may be removed
from this representation, while for $k > \frac{m}{2}$, the $o$-term dominates.
Moreover, if $\frac{m+2}{2} < k \leq m$, we are left with this term, only,
that is,
\[
D(Z_n) = o \biggl(\frac{1}{n^{(m-2)/2} (\log n)^{(m-d)/2} } \biggr).
\]

When $k = 3$, there is no restriction on the cumulants in Corollary \ref
{Corollary6.2},
and (\ref{equ6.11}) becomes
\[
D(Z_n) = \frac{1}{2n} \sum_{|\nu| = 3}
\frac{\gamma_\nu^2}{\nu!} + O \biggl(\frac{1}{n^2} \biggr) + o
\biggl(
\frac{1}{n^{(m-2)/2} (\log
n)^{(m-d)/2} } \biggr).
\]

If $\E|X_1|^4 < +\infty$, we get $D(Z_n) = O(1/n)$ for $d \leq4$,
and the weaker bound $D(Z_n) = o((\log n)^{(d-4)/2}/n)$ for $d \geq5$.
However, if $\E|X_1|^5 < +\infty$, we always have $D(Z_n) = O(1/n)$
regardless of the dimension $d$.

Technically, this slight difference between conclusions for different dimensions
is due to the dimension-dependent asymptotic
$\int_{|x| > T} |x|^2 \varphi(x) \,dx \sim C_d T^d e^{-T^2/2}$.

\begin{Remark*} In case of discrete distributions when $X_1$ takes integer
values, asymptotics for $D(S_n)$ were studied by Vilenkin and D'yachkov
\cite{V-D},
who used an Edgeworth-type expansion for probabilities $\P\{S_n = k\}$
in the corresponding local limit theorem.
\end{Remark*}

\section{Convolutions of mixtures of normal laws}\label{sec7}

Is the asymptotic description of $D(Z_n)$ in Theorem~\ref{Theorem1.1}
still optimal,
if no
expansion terms of order $n^{-j}$ are present?
This is exactly the case for $2 \leq s < 4$.

In order to answer the question, we examine a special class of probability
distributions that can be described as mixtures of normal laws on the
real line
with mean zero. They have densities of the form
%
%
\begin{equation}\label{equ7.1}
p(x) = \int_0^{+\infty} \varphi_\sigma(x)
\,dP(\sigma)\qquad (x \in\R),
\end{equation}
where $P$ is a (mixing) probability measure on the positive half-axis
$(0,+\infty)$,
and where
\[
\varphi_\sigma(x) = \frac{1}{\sigma\sqrt{2\pi}} e^{-x^2/(2\sigma^2)}
\]
is the density of the normal law with mean zero and variance $\sigma^2$
[as usual, we write $\varphi(x)$ in the standard normal case with
$\sigma=1$].

Equivalently, let $p(x)$ denote the density of the random variable $X_1
= \rho Z$,
where the factors $Z \sim N(0,1)$ and $\rho>0$ (with the distribution
$P$) are independent. Such distributions appear naturally, for example,
as limit laws of sums with randomized length; cf., for example,~\cite{B-G}.

For densities such as (\ref{equ7.1}), we need a refinement of the local
limit theorem
for convolutions, described in the expansions (\ref{equ3.5}) and (\ref
{equ3.6}). More precisely,
our aim is to find a representation with an essentially smaller
remainder term
compared to $o(n^{-(s-2)/2})$.

Thus, let $X_1,X_2,\ldots$ be independent random variables, having a
common density
$p(x)$ as in (\ref{equ7.1}), and let $p_n(x)$ denote the density of the
normalized sum
$Z_n = (X_1 + \cdots+ X_n)/\sqrt{n}$.
If $X_1 = \rho Z$, where $Z \sim N(0,1)$ and $\rho>0$ are independent, then
$\E X_1^2 = \E\rho^2$ and more generally,
\[
\E|X_1|^s = \beta_s \E\rho^s =
\beta_s \int_0^{+\infty}
\sigma^s \,dP(\sigma),
\]
where $\beta_s$ denotes the $s$th absolute moment of $Z$.

Note that $p(x)$ is unimodal with mode at the origin, and
$p(0) = \E\frac{1}{\rho\sqrt{2\pi}}$.
If $\rho\geq\sigma_0 > 0$, the density is bounded, and therefore the entropy
$h(X_1)$ is finite.

%
\begin{Proposition}\label{Proposition7.1} Assume that $\E\rho^2 = 1$,
$\E\rho^s <
+\infty$
$(2 < s \leq4)$. If $\P\{\rho\geq\sigma_0\} = 1$ with some constant
$\sigma_0 > 0$, then uniformly over all $x$,
%
%
\begin{equation}\label{equ7.2}
p_n(x) = \varphi(x) + n \int_0^{+\infty}
\bigl(\varphi_{\sigma_n}(x) - \varphi(x) \bigr) \,dP(\sigma) + O
\biggl(
\frac{1}{n^{s-2}} \biggr),
\end{equation}
where $\sigma_n = \sqrt{1 + \frac{\sigma^2 - 1}{n}}$.
\end{Proposition}

Of course, when $\E\rho^s < +\infty$ for $s > 4$, the proposition may
be still
applied, but with $s=4$. In this case (\ref{equ7.2}) has a remainder
term of order
$O(\frac{1}{n^2})$. Note that necessarily $\sigma_0 \leq1$ under
the condition $\E\rho^2 = 1$.

The function $p_n$ may also be described as the density of
$Z_n = \sqrt{\frac{\rho_1^2 + \cdots+ \rho_n^2}{n}} Z$, where $\rho
_k$ are
independent copies of $\rho$ (independent of $Z$ as well). This represention
already indicates the closeness of $p_n$ and $\varphi$ and suggests to appeal
to the law of large numbers. However, we shall choose
a different approach based on the characteristic functions of $Z_n$.

Obviously, the characteristic function of $X_1$ is given by
\[
v(t) = \E e^{itX_1} = \E e^{-\rho^2 t^2/2}\qquad (t \in\R).
\]
Using Jensen's inequality and the assumption $\rho\geq\sigma_0 > 0$,
we get
a two-sided estimate
%
%
\begin{equation}\label{equ7.3}
e^{-t^2/2} \leq v(t) \leq e^{-\sigma_0^2 t^2/2}.
\end{equation}
In particular, the function $\psi(t) = e^{t^2/2} v(t) - 1$
is nonnegative for all $t$ real.

%
\begin{Lemma}\label{Lemma7.2} If $\E\rho^2 = 1$, $M_s = \E\rho^s <
+\infty$
$(2 \leq s \leq4)$, then for all \mbox{$|t| \leq1$},
\[
0 \leq\psi(t) \leq M_s |t|^s.
\]
\end{Lemma}

\begin{pf} We may assume $0 < t \leq1$.
Write $\psi(t) = \E(e^{-(\rho^2 - 1) t^2/2} - 1 )$.
The expression under the expectation sign is nonpositive for $\rho t >
1$, hence
\[
\psi(t) \leq\E\bigl(e^{-(\rho^2 - 1) t^2/2} - 1 \bigr) 1_{\{\rho
\leq
1/t\}}.\vadjust{\goodbreak}
\]
Let $x = -(\rho^2 - 1) t^2$. Clearly,
$|x| \leq1$ for $\rho\leq1/t$. Using $e^x \leq1 + x + x^2$ ($|x|
\leq1$)
and $\E\rho^2 = 1$, we get
%
%
\begin{eqnarray}\label{equ7.4}
\psi(t) & \leq& - \frac{t^2}{2} \E\bigl(\rho^2 - 1\bigr)
1_{\{\rho\leq1/t\}} + \frac{t^4}{4} \E\bigl(\rho^2 - 1
\bigr)^2 1_{\{\rho\leq1/t\}}
\nonumber\\[-8pt]\\[-8pt]
& = & \frac{t^2}{2} \E\bigl(\rho^2 - 1\bigr) 1_{\{\rho> 1/t\}}
+ \frac{t^4}{4} \E\bigl(\rho^2 - 1\bigr)^2
1_{\{\rho\leq1/t\}}.\nonumber
\end{eqnarray}
The last expectation is equal to
\begin{eqnarray*}
&&
\E\rho^4 1_{\{\rho\leq1/t\}} + 2 \E\bigl(\rho^2 - 1\bigr)
1_{\{\rho> 1/t\}} - \P\{\rho\leq1/t\} \\
&&\qquad \leq \E\rho^4
1_{\{\rho\leq1/t\}} + 2 \E\rho^2 1_{\{\rho> 1/t\}} - 1
\\
&&\qquad \leq \E\rho^4 1_{\{\rho\leq1/t\}} + \E\rho^2
1_{\{\rho> 1/t\}}.
\end{eqnarray*}
Together with (\ref{equ7.4}), this gives
%
%
\begin{equation}\label{equ7.5}
\psi(t) \leq\frac{3t^2}{4} \E\rho^2 1_{\{\rho> 1/t\}} +
\frac{t^4}{4} \E\rho^4 1_{\{\rho\leq1/t\}}.
\end{equation}
Finally,\vspace*{1pt}
$\E\rho^2 1_{\{\rho> 1/t\}} \leq\E\rho^s t^{s-2} 1_{\{\rho> 1/t\}}
\leq M_s t^{s-2}$ and
$\E\rho^4 1_{\{\rho\leq1/t\}} \leq
\E\rho^s t^{s-4} \*1_{\{\rho\leq1/t\}} \leq M_s t^{s-4}$.
It remains to use these estimates in (\ref{equ7.5}), and Lemma~\ref
{Lemma7.2} is proved.
\end{pf}

\begin{pf*}{Proof of Proposition~\ref{Proposition7.1}}
The characteristic functions $v_n(t) = v(\frac{t}{\sqrt{n}})^n$ of
$Z_n$ are real-valued and admit, by (\ref{equ7.3}), similar bounds
%
%
\begin{equation}\label{equ7.6}
e^{-t^2/2} \leq v_n(t) \leq e^{-\sigma_0^2 t^2/2}.
\end{equation}
In particular, one may apply the inverse Fourier transform to represent
the density
of $Z_n$ as
\[
p_n(x) = \frac{1}{2\pi} \int_{-\infty}^{+\infty}
e^{-itx} v_n(t) \,dt = \frac{1}{2\pi} \int
_{-\infty}^{+\infty} e^{-itx - t^2/2} \bigl(1 + \psi(t/
\sqrt{n}) \bigr)^n \,dt.
\]
Letting $T_n = \frac{4}{\sigma_0} \log n$, we split the integral into
the two regions,
defined by
\[
I_1 = \int_{|t| \leq T_n} e^{-itx}
v_n(t) \,dt,\qquad I_2 = \int_{|t| > T_n}
e^{-itx} v_n(t) \,dt.
\]
By the upper bound in (\ref{equ7.6}),
%
%
\begin{equation}\label{equ7.7}
|I_2| \leq\int_{|t| > T_n} e^{-\sigma_0^2 t^2/2} \,dt \leq
\frac{\sqrt{2\pi}}{\sigma_0} e^{-\sigma_0^2 T_n^2/2} = \frac
{\sqrt{2\pi}}{\sigma_0 n^8}.
\end{equation}
In the interval $|t| \leq T_n$, by Lemma~\ref{Lemma7.2},
$\psi(\frac{t}{\sqrt{n}}) \leq\frac{M_s |t|^s}{n^{s/2}} \leq\frac{1}{n}$,
for all $n \geq n_0$. But for $0 \leq\ep\leq\frac{1}{n}$, there is
the simple estimate $0 \leq(1 + \ep)^n - 1 - n \ep\leq2 (n\ep)^2$.
Hence, once more by Lemma~\ref{Lemma7.2},
\begin{eqnarray*}
0 &\leq&\bigl(1 + \psi(t/\sqrt{n}) \bigr)^n - 1 - n \psi(t/\sqrt{n})\\
&\leq&2 \bigl(n \psi(t/\sqrt{n}) \bigr)^2 \leq2M_s^2
\frac{|t|^{2s}}{n^{s-2}}\qquad (n \geq n_0).
\end{eqnarray*}
This gives
%
%
\begin{equation}\label{equ7.8}\qquad
\biggl|I_1 - \int_{|t| \leq T_n} e^{-itx - t^2/2} \bigl(1 +
n\psi(t/\sqrt{n}) \bigr) \,dt \biggr| \leq\frac{2M_s^2}{n^{s-2}} \int
_{-\infty}^{+\infty}
|t|^{2s} e^{-t^2/2} \,dt.
\end{equation}
In addition,
\begin{eqnarray*}
&&
\biggl|\int_{|t| > T_n} e^{-itx - t^2/2} \bigl(1 + n\psi(t/\sqrt{n})
\bigr) \,dt \biggr| \\
&&\qquad\leq\int_{|t| > T_n} e^{-t^2/2} \,dt + n \int
_{|t| > T_n} e^{- t^2/2} \psi(t/\sqrt{n}) \,dt.
\end{eqnarray*}
Here, the first integral on the right-hand side is of order $O(n^{-8})$.
To estimate the second one, recall that, by (\ref{equ7.3}),
$\psi(t) = e^{t^2/2} v(t) - 1 \leq e^{(1 - \sigma_0^2) t^2/2}$. Hence,
$\psi(t/\sqrt{n}) \leq e^{(1 - \sigma_0^2) t^2/2}$ and
\[
\int_{|t| > T_n} e^{-t^2/2} \psi(t/\sqrt{n}) \,dt \leq\int
_{|t| > T_n} e^{-\sigma_0^2 t^2/2} \,dt \leq\frac{\sqrt{2\pi
}}{\sigma_0 n^8}.
\]
Together with (\ref{equ7.7}) and (\ref{equ7.8}) these bounds imply that
\[
p_n(x) = \frac{1}{2\pi} \int_{-\infty}^{+\infty}
e^{-itx - t^2/2} \bigl(1 + n\psi(t/\sqrt{n}) \bigr) \,dt + O \biggl(
\frac{1}{n^{s-2}} \biggr)
\]
uniformly over all $x$. It remains to note that
\begin{eqnarray*}
\frac{1}{2\pi} \int_{-\infty}^{+\infty}
e^{-itx - t^2/2} \psi(t/\sqrt{n}) \,dt & = & \frac{1}{2\pi} \int
_{-\infty}^{+\infty} e^{-itx - t^2/2} \bigl(e^{t^2/2n}
v(t/\sqrt{n}) - 1 \bigr) \,dt
\\
& = & \int_0^{+\infty} \bigl(\varphi_{\sigma_n}(x)
- \varphi(x) \bigr) \,dP(\sigma).
\end{eqnarray*}
Proposition~\ref{Proposition7.1} is proved.
\end{pf*}

%
\begin{Remark}\label{Remark7.3}
An inspection of (\ref{equ7.5}) shows that, in the case $2 < s < 4$,
Lemma~\ref{Lemma7.2}
may slightly be sharpened to $\psi(t) = o(|t|^s)$. Correspondingly, the
$O$-relation in Proposition~\ref{Proposition7.1} can be replaced with
an $o$-relation.
This improvement is convenient, but not crucial for the proof of
Theorem~\ref{Theorem1.3}.
\end{Remark}

\section{\texorpdfstring{Lower bounds. Proof of Theorem \protect\ref{Theorem1.3}}
{Lower bounds. Proof of Theorem 1.3}}\label{sec8}

Let $X_1,X_2,\ldots$ be independent random variables with a common
density of the
form
\[
p(x) = \int_0^{+\infty} \varphi_\sigma(x)
\,dP(\sigma),\qquad x \in\R.\vadjust{\goodbreak}
\]
Equivalently, let $X_1 = \rho Z$ with independent random variables $Z
\sim N(0,1)$
and $\rho>0$ having distribution $P$.

A basic tool for proving Theorem~\ref{Theorem1.3} will be the following
lower bound
on the
entropic distance to Gaussianity for the partial sums $S_n = X_1 +
\cdots+ X_n$.

%
\begin{Proposition}\label{Proposition8.1} Let $\E\rho^2 = 1$, $\E\rho^s
< +\infty$
$(2 < s < 4)$ and $\P\{\rho\geq\sigma_0\} = 1$ with $\sigma_0 > 0$.
Assume that, for some $\gamma> 0$,
%
%
\begin{equation}\label{equ8.1}
\liminf_{n \rightarrow\infty} n^{s - {1/2}} \int_{n^{
{1/2} + \gamma}}^{+\infty}
\frac{1}{\sigma} \,dP(\sigma) > 0.
\end{equation}
Then with some absolute constant $c>0$ and some constant $\delta>0$,
%
%
\begin{equation}\label{equ8.2}
D(S_n) \geq c n \log n \P\{\rho\geq\sqrt{n\log n} \} + o \biggl(
\frac{1}{n^{({s-2})/{2} + \delta}} \biggr).
\end{equation}
\end{Proposition}

In fact, in (\ref{equ8.2}) one may take any positive number
$\delta< \min\{\gamma s, \frac{s-2}{2}\}$.

\begin{pf*}{Proof of Proposition~\ref{Proposition8.1}}
By Proposition~\ref{Proposition7.1} and Remark~\ref{Remark7.3},
uniformly over all $x$,
%
%
\begin{equation}\label{equ8.3}
p_n(x) = \varphi(x) + n \int_0^{+\infty}
\bigl(\varphi_{\sigma_n}(x) - \varphi(x) \bigr) \,dP(\sigma) + o
\biggl(
\frac{1}{n^{s-2}} \biggr),
\end{equation}
where $p_n$ is the density of $S_n/\sqrt{n}$ and
$\sigma_n = \sqrt{1 + \frac{\sigma^2 - 1}{n}}$.

Define the sequence
\[
N_n = \frac{n^{{1}/{2} + \gamma}}{5 \sqrt{\log n}}
\]
for $n$ large enough (so that $N_n \geq1$). By Chebyshev's inequality,
%
%
\begin{equation}\label{equ8.4}
\P\{\rho\geq N_n\} \leq5^s M_s
\frac{\log^2 n}{n^{({1}/{2} +
\gamma) s}} = o \biggl(\frac{1}{n^{{s/2} + \delta}} \biggr),
\qquad 0 < \delta< \gamma s.
\end{equation}

Using $u \log u \geq u-1$ ($u \geq0$) and applying (\ref{equ8.3}), we
may write
%
%
\begin{eqnarray}\label{equ8.5}
I_n &\equiv& \int_{|x| \leq4 \sqrt{\log n}}
p_n(x) \log\frac{p_n(x)}{\varphi(x)} \,dx \nonumber\\
& \geq& \int_{|x| \leq4
\sqrt{\log n}}
\bigl(p_n(x) - \varphi(x) \bigr) \,dx
\\
& \geq & n \int_0^{+\infty} \int
_{|x| \leq4 \sqrt{\log n}} \bigl(\varphi_{\sigma_n}(x) - \varphi
(x) \bigr) \,dx
\,dP(\sigma) - C \frac{\sqrt{\log n}}{n^{s-2}}\nonumber
\end{eqnarray}
with some constant $C$.

Note that $\sigma_n < 1$ for $\sigma< 1$, and thus, for any $T>0$,
\[
\int_{|x| \leq T} \bigl(\varphi_{\sigma_n}(x) - \varphi(x)
\bigr) \,dx = 2 \bigl(\Phi(T/\sigma_n) - \Phi(T) \bigr) > 0,
\]
where $\Phi$ denotes the distribution function of the standard normal law.
Hence, the outer integral in (\ref{equ8.5}) may be restricted to the range
$\sigma\geq1$. Moreover, by (\ref{equ8.4}), one may also restrict this
integral, even
to the range $\sigma\geq N_n$. More precisely, (\ref{equ8.4}) gives
\[
n \biggl|\int_{N_n}^{+\infty} \int_{|x| \leq4 \sqrt{\log n}}
\bigl(\varphi_{\sigma_n}(x) - \varphi(x) \bigr) \,dx \,dP(\sigma) \biggr|
\leq n \P\{
\rho\geq N_n\} = o \biggl(\frac{1}{n^{({s-2})/{2} + \delta}}
\biggr).
\]
Comparing this relation with (\ref{equ8.5}) and imposing the additional
requirement
$\delta< \frac{s-2}{2}$, we get
%
%
\begin{eqnarray}\label{equ8.6}\quad
I_n & \geq& n \int_1^{N_n} \int
_{|x| \leq4 \sqrt{\log n}} \bigl(\varphi_{\sigma_n}(x) - \varphi
(x) \bigr) \,dx
\,dP(\sigma) + o \biggl(\frac{1}{n^{({s-2})/{2} + \delta}} \biggr)
\nonumber\\[-8pt]\\[-8pt]
& = & -2n \int_1^{N_n} \int
_{{4}\sqrt{\log n}/{\sigma_n}}^{4 \sqrt{\log n}} \varphi(x)
\,dx \,dP(\sigma) + o \biggl(
\frac{1}{n^{({s-2})/{2} + \delta
}} \biggr).\nonumber
\end{eqnarray}

Now, let us estimate $p_n(x)$ from below in the region
$4 \sqrt{\log n} \leq|x| \leq n^{\gamma}$. If $|x| \geq4 \sqrt
{\log n}$,
it follows from (\ref{equ8.3}) that
%
%
\begin{equation}\label{equ8.7}
p_n(x) = n \int_0^{+\infty}
\varphi_{\sigma_n}(x) \,dP(\sigma) + o \biggl(\frac{1}{n^{s-2}}
\biggr).
\end{equation}

Consider the function
\[
g_n(x) = \int_0^{+\infty}
\frac{\varphi_{\sigma_n}(x)}{\varphi(x)} \,dP(\sigma).
\]
Note that $1 \leq\sigma_n \leq\sigma$ for $\sigma\geq1$. In this case,
the ratio $\frac{\varphi_{\sigma_n}(x)}{\varphi(x)}$ is nonincreasing
in $x \geq0$.
Moreover, for $\sigma\geq\sqrt{3n + 1}$, we have
$\sigma_n^2 = 1 + \frac{\sigma^2 - 1}{n} \geq4$, so
$1 - \frac{1}{\sigma_n^2} \geq\frac{3}{4}$. Hence, for $|x| \geq4
\sqrt{\log n}$,
\[
\frac{\varphi_{\sigma_n}(x)}{\varphi(x)} = \frac{1}{\sigma_n}
e^{{x^2}(1 - {1}/{\sigma_n^2})/{2}} \geq
\frac{n^6}{\sigma}.
\]
Therefore,
\[
g_n(x) \geq n^6 \int_{\sqrt{3n+1}}^{+\infty}
\frac{1}{\sigma} \,dP(\sigma).
\]
But by assumption (\ref{equ8.1}), the last expression tends to infinity
with $n$,
so for all $n$ large enough, $g_n(x) \geq2$ in the interval $|x| \geq
4 \sqrt{\log n}$.

Furthermore, if $\sigma\geq|x|\sqrt{n}$, then
$\sigma_n^2 = 1 + \frac{\sigma^2 - 1}{n} \geq x^2$,
so $\frac{x^2}{2 \sigma_n^2} \leq\frac{1}{2}$. On the other hand,
\[
\sigma_n^2 < 1 + \frac{\sigma^2}{n} =
\frac{n + \sigma^2}{n} \leq\frac{{\sigma^2}/{x^2} + \sigma
^2}{n} \leq\frac{2\sigma^2}{n},
\]
since
$|x| \geq4\sqrt{\log n} > 1$ for $n \geq2$. The two estimates give
\[
\varphi_{\sigma_n}(x) = \frac{1}{\sigma_n \sqrt{2\pi}}
e^{-x^2/2\sigma_n^2} \geq
\frac{\sqrt{n}}{6\sigma}.
\]
Therefore, whenever $4 \sqrt{\log n} \leq|x| \leq n^\gamma$,
\[
n \int_0^{+\infty} \varphi_{\sigma_n}(x) \,dP(
\sigma) \geq\frac{n^{3/2}}{6} \int_{|x|\sqrt{n}}^{+\infty}
\frac{1}{\sigma} \,dP(\sigma) \geq\frac{n^{3/2}}{6} \int_{n^{
{1}/{2} + \gamma}}^{+\infty}
\frac
{1}{\sigma} \,dP(\sigma).
\]
By assumption (\ref{equ8.1}), the last expression and therefore the left
integral are larger than $\frac{c}{n^{s-2}}$ with some constant $c>0$.
Consequently, the remainder term in (\ref{equ8.7}) is indeed smaller,
so that
for all $n$ large enough, we may write, for example,
\[
p_n(x) \geq0.8 n \int_0^{+\infty}
\varphi_{\sigma_n}(x) \,dP(\sigma) = 0.8 n g_n(x) \varphi(x)
\qquad\bigl(4 \sqrt{\log n} \leq|x| \leq n^\gamma\bigr).
\]

Since $g_n(x) \geq2$ for $|x| \geq4 \sqrt{\log n}$ with large $n$,
we have in this region $\frac{p_n(x)}{\varphi(x)} \geq1.6 n > n$, thus
\[
p_n(x) \log\frac{p_n(x)}{\varphi(x)} \geq p_n(x) \log n \geq
0.8 n \log n \int_0^{+\infty} \varphi_{\sigma_n}(x)
\,dx \,dP(\sigma).
\]
Hence,
%
%
\begin{eqnarray}\label{equ8.8}
&&
\int_{4 \sqrt{\log n} \leq|x| \leq n^\gamma} p_n(x) \log\frac
{p_n(x)}{\varphi(x)} \,dx \nonumber\\
&&\qquad
\geq 0.8 n \log n \int_0^{+\infty} \int
_{4 \sqrt{\log n} \leq|x| \leq n^\gamma}\varphi_{\sigma_n}(x) \,dx
\,dP(\sigma)
\\
&&\qquad = 1.6 n \log n \int_0^{+\infty}
\int_{{4}\sqrt{\log n}/{\sigma_n}}^{{n^\gamma}/{\sigma
_n}} \varphi(x) \,dx \,dP(\sigma).
\nonumber
\end{eqnarray}
At this point, it is useful to note that
$\frac{n^\gamma}{\sigma_n} \geq4 \sqrt{\log n}$,
as long as $\sigma\leq N_n$ with $n$ large enough. Indeed, in this case
$\sigma_n^2 \leq(1 - \frac{1}{n}) + \frac{N_n^2}{n} < 1 + \frac
{n^{2\gamma}}{25 \log n}$,
so
\[
(4\sigma_n \sqrt{\log n} )^2 \leq16 \log n \biggl(1 +
\frac{n^{2\gamma}}{25 \log n} \biggr) < n^{2\gamma}
\]
for all $n$ large enough. Hence, from (\ref{equ8.8}),
\[
\int_{4 \sqrt{\log n} \leq|x| \leq n^\gamma} p_n(x) \log\frac
{p_n(x)}{\varphi(x)} \,dx
\geq1.6 n \log n \int_0^{N_n} \int
_{{4}\sqrt{\log n}/{\sigma_n}}^{4\sqrt{\log n}} \varphi(x) \,dx
\,dP(\sigma).\vadjust{\goodbreak}
\]

But the last expression dominates the double integral in (\ref{equ8.6})
with a factor of~$2n$. Therefore, combining the above estimate with
(\ref{equ8.6}), we get
\begin{eqnarray*}
\int_{|x| \leq n^\gamma} p_n(x) \log\frac{p_n(x)}{\varphi(x)} \,dx
&\geq&
1.4 n \log n \int_0^{N_n} \int
_{{4}\sqrt{\log n}/{\sigma_n}}^{4\sqrt{\log n}} \varphi(x) \,dx
\,dP(\sigma) \\
&&{}+ o \biggl(
\frac{1}{n^{({s-2})/{2} + \delta}} \biggr).
\end{eqnarray*}

Finally, we may extend the outer integral on the right-hand side
to all values $\sigma> 0$ by noting that, by (\ref{equ8.4}),
\[
n \log n \int_{N_n}^{+\infty} \int_{{4}\sqrt{\log
n}/{\sigma_n}}^{4\sqrt{\log n}}
\varphi(x) \,dx \,dP(\sigma) \leq n \log n \P\{\rho> N_n\} = o \biggl(
\frac{1}{n^{({s-2})/{2} + \delta}} \biggr).
\]
Hence,
%
%
\begin{eqnarray}\label{equ8.9}\qquad
\int_{|x| \leq n^\gamma} p_n(x) \log\frac{p_n(x)}{\varphi(x)} \,dx
&\geq&
1.4 n \log n \int_0^{+\infty} \int
_{{4}\sqrt{\log n}/{\sigma_n}}^{4\sqrt{\log n}} \varphi(x) \,dx
\,dP(\sigma)\nonumber\\[-8pt]\\[-8pt]
&&{} + o \biggl(
\frac{1}{n^{({s-2})/{2} + \delta}} \biggr).\nonumber
\end{eqnarray}

For the remaining values $|x| \geq n^\gamma$, one can just use the property
$u \log u \geq- \frac{1}{e}$ to get a simple lower bound
\begin{eqnarray*}
\int_{|x| > n^\gamma} p_n(x) \log\frac{p_n(x)}{\varphi(x)} \,dx &
\geq& \int_{|x| > n^\gamma, p_n(x) \leq\varphi(x)} p_n(x) \log
\frac
{p_n(x)}{\varphi(x)} \,dx
\\
& \geq& - \frac{1}{e} \int_{|x| > n^\gamma, p_n(x) \leq\varphi
(x)} \varphi(x) \,dx \geq
- e^{-n^{2\gamma}/2}.
\end{eqnarray*}
Together with (\ref{equ8.9}) this yields
\begin{eqnarray*}
\int_{-\infty}^{+\infty} p_n(x) \log
\frac{p_n(x)}{\varphi(x)} \,dx &\geq& 1.4 n \log n \int_0^{+\infty}
\int_{{4}\sqrt{\log n}/{\sigma_n}}^{4\sqrt{\log n}} \varphi
(x) \,dx \,dP(\sigma) \\
&&{}
+ o \biggl(\frac{1}{n^{({s-2})/{2} + \delta}} \biggr).
\end{eqnarray*}

To simplify, finally note that $\frac{4}{\sigma_n} \sqrt{\log n}
\leq4$
for $\sigma\geq\sqrt{n\log n}$. In this case the last integral is separated
from zero (for large $n$), hence with some absolute constant $c>0$
\[
\int_{-\infty}^{+\infty} p_n(x) \log
\frac{p_n(x)}{\varphi(x)} \,dx \geq c n \log n \P\{\rho\geq\sqrt
{n\log n} \} + o \biggl(
\frac{1}{n^{({s-2})/{2} + \delta}} \biggr).
\]
This is exactly the required inequality (\ref{equ8.2}) and Proposition
\ref{Proposition8.1} is proved.
\end{pf*}

\begin{pf*}{Proof of Theorem~\ref{Theorem1.3}}
Given $\eta> 0$, one may apply Proposition~\ref{Proposition8.1} to the
probability
measure $P$ with density
\[
\frac{dP(\sigma)}{d\sigma} = \frac{c}{\sigma^{s+1} (\log\sigma
)^\eta},\qquad \sigma> 2,
\]
and extending it to an interval $[\sigma_0,2]$ to meet the requirement
$\int_{\sigma_0}^{+\infty} \sigma^2 \,dP(\sigma) = 1$ (with some
$0 < \sigma_0 < 1$ and a positive normalizing constant $c = c_{\eta,s}$).
It is easy to see that in this case condition (\ref{equ8.1}) is
fulfilled for
$0 < \gamma< \frac{s-2}{2(s+1)}$. In addition, if $\rho$ has
the distribution $P$, we have
\[
\P\{\rho\geq\sigma\} \geq \operatorname{const} \frac{1}{\sigma^s (\log
\sigma
)^\eta}
\]
for all $\sigma$ large enough.
Hence, by taking $\sigma= \sqrt{n \log n}$, (\ref{equ8.2}) provides the desired
lower bound.
\end{pf*}

\begin{Remark*} In case $s=2$ (i.e., with minimal moment assumptions),
the mixtures of the normal laws with discrete mixing measures $P$
were used by Matskyavichyus~\cite{M} in the central limit theorem in terms
of the
Kolmogorov distance. Namely, it is shown that, for any prescribed sequence
$\ep_n \rightarrow0$, one may choose $P$ such that
$\Delta_n = \sup_x |F_n(x) - \Phi(x)| \geq\ep_n$ for all $n$ large enough
(where $F_n$ is the distribution function of $Z_n$). In view of the Pinsker-type
inequality, one may conclude that
\[
D(Z_n) \geq\tfrac{1}{2} \Delta_n^2 \geq
\tfrac{1}{2} \ep_n^2.
\]
Therefore, $D(Z_n)$ may decay at an arbitrarily slow rate.
\end{Remark*}



\printaddresses

\end{document}